\documentclass[journal,draftclsnofoot,onecolumn,11pt]{IEEEtran}
\usepackage[usenames,dvipsnames]{color}
\usepackage{graphicx}
\DeclareGraphicsExtensions{.pdf,.jpeg,.png,.eps}
\graphicspath{{./figures/}}

\usepackage[cmex10]{amsmath}
\usepackage{amsfonts}
\usepackage{amsthm}

\usepackage{amssymb}
\usepackage{mathrsfs}    

\newcommand{\be}{\begin{equation}}
\newcommand{\ee}{\end{equation}}
\renewcommand{\Vec}[1]{\boldsymbol{#1}}
\newcommand{\Mat}[1]{\boldsymbol{#1}}

\renewcommand{\Re}{\mathop{\mathrm{Re}}}

\newtheorem{theorem}{Theorem}
\newtheorem{cor}{Corollary}
\newtheorem{prop}{Proposition}

\newcounter{saveeqn}
\newcommand{\alpheqn}{\stepcounter{equation}%
\setcounter{saveeqn}{\value{equation}}\setcounter{equation}{0}%
\renewcommand{\theequation}
            {\arabic{saveeqn}\alph{equation}}}
\newcommand{\reseteqn}{\setcounter{equation}{\value{saveeqn}}%
\renewcommand{\theequation}{\arabic{equation}}}

\newcommand{\oldtext}[1]{\noindent [old text is deleted]\\}

\newcommand{\delete}[1]{}

\usepackage{ulem}

\DeclareMathOperator*{\argmin}{\arg\!\min}

\usepackage{dsfont}

\usepackage{psfrag}		
\usepackage{epstopdf}	
\usepackage{NumericPlots}

\begin{document}

\title{Using the LASSO's Dual for Regularization in Sparse Signal Reconstruction from Array Data}

\author{
Christoph F. Mecklenbr\"auker, 
Peter Gerstoft, 
Erich Z\"ochmann
\thanks{\today}
\thanks{submitted to IEEE Transactions on Signal Processing on 09-Aug-2015, Manuscript ID: T-SP-19417-2015}
\thanks{Christoph F. Mecklenbr\"auker and Erich Z\"ochmann are with Institute of Telecommunications,
Vienna University of Technology
1040 Vienna, Austria, {\tt cfm@ieee.org}}
\thanks{Peter Gerstoft is with University of California San Diego, La Jolla, CA 92093-0238, USA
http://www.mpl.ucsd.edu/people/pgerstoft}
\thanks{\copyright ~ IEEE 2015}
}

\maketitle

\begin{abstract}
Waves from a sparse set of source hidden in additive noise are observed by a sensor array.
We treat the estimation of the sparse set of sources as a generalized complex-valued LASSO problem.
The corresponding dual problem is formulated and it is shown that the dual solution is useful for selecting the regularization parameter of the LASSO when the number of sources is given. 
The solution path of the complex-valued LASSO is analyzed. 
For a given number of sources,  the corresponding regularization parameter is determined by an order-recursive algorithm and two iterative algorithms that are based on a further approximation.
Using this regularization parameter, the DOAs of all sources are estimated.
\end{abstract}
{\IEEEkeywords sparsity, generalized LASSO, duality theory}
\IEEEpeerreviewmaketitle

\section{Introduction}

This paper contributes to the area of sparse signal estimation for sensor array processing. Sparse signal estimation techniques retrieve a signal vector from an undercomplete set of noisy measurements when the signal vector is assumed to have only few nonzero components at unknown positions.
Research in this area was spawned by the Least Absolute Shrinkage and Selection Operator (LASSO) \cite{Tibshirani1996}.
In the related field of compressed sensing, this sparse signal reconstruction problem is known as the atomic decomposition problem \cite{ChenDonohoSaunders1998}.
The early results for sparse signals \cite{gorodnitsky1997,candes06,donoho06}
have been extended to compressible (approximately sparse) signals and sparse signals buried in noise \cite{Fuchs2005,malioutov05,Donoho-Elad-Temlyakov2006,Tropp2006,Elad2010} which renders the framework applicable to problems in array processing.  

Similar to \cite{Tibshirani2011,Andrade2010}, the LASSO is generalized and formulated for complex-valued observations acquired from a sensor array. 
It is shown here that the corresponding  dual vector is interpretable as the output of a weighted matched filter (WMF) acting on the residuals of the linear observation model, cf. \cite{Donoho2012}.

The regularization parameter $\mu$ in LASSO defines the trade-off between the model fit and the estimated sparsity order $K$ given by the number of estimated nonzero signal components.
When the sparsity order $K_0$ is given, choosing a suitable value for the LASSO regularization parameter $\mu$ remains a challenging task.
The homotopy techniques \cite{OsbornePresnellTurlach2000,Panahi-Viberg2012,Koldovsky2014} provide an approach to sweep over a range of $\mu$ values to select the signal estimate with the given $K_0$.

The  maximum magnitudes of the dual vector can be used for selecting the regularization parameter of the generalized LASSO. This is the basis for an order-recursive algorithm to solve the sparse signal reconstruction problem \cite{Koldovsky2014,Panahi-Viberg2011,MGPV2013} for the given $K_0$.
In this work,  a fast and efficient choice of $\mu$ is proposed for direction of arrival estimation from array data. The choice exploits the sidelobe levels of the array's beampattern. We motivate this choice after proving several  relations between the regularization parameter $\mu$, the LASSO residuals, and the LASSO's dual solution.

The main achievements of this work are summarized as follows: We extend the convex duality theory \cite{Tibshirani2011} from the real-valued to the complex-valued case and formulate the corresponding dual problem to the complex-valued LASSO. We show that the dual solution is useful for selecting the regularization parameter of the LASSO. Three signal processing algorithms are formulated and evaluated to support our theoretical results and claims.

\subsection{Notation}
Matrices $\Mat{A},\Mat{B},\ldots$ and vectors $\Vec{a},\Vec{b},\ldots$ are complex-valued and denoted by boldface letters. The zero vector is $\Vec{0}$. The Hermitian transpose, inverse, and Moore-Penrose pseudo inverse are denoted as $\Mat{X}^H,\Mat{X}^{-1},\Mat{X}^+$ respectively. We abbreviate $\Mat{X}^{-H}=(\Mat{X}^H)^{-1}$. The complex vector space of dimension $N$ is written as $\mathbb{C}^N$. $\mathcal{N}(\Mat{A})$ is the null space of $\Mat{A}$ and $\mathop{\mathrm{span}}(\Mat{A})$ denotes the linear hull of $\Mat{A}$. The projection onto $\mathop{\mathrm{span}}(\Mat{A})$ is $\Mat{P}_{\!_{\!\Mat{A}}}$. The $\ell_p$-norm is written as $\|\cdot\|_p$. For a vector $\Vec{x}\in\mathbb{C}^M$, $\|\Vec{x}\|_{\infty} = \max\limits_{1\le m\le M} |x_m|$, for a matrix $\Mat{X}\in\mathbb{C}^{N\times M}$, we define $\|\Mat{X}\|_{\infty} = \max\limits_{1\le n\le N} \max\limits_{1\le m\le M} |X_{nm}|$.

\section{Problem formulation}

We start from the following problem formulation: Let $\Vec{y}\in\mathbb{C}^N$ and $\Mat{A}\in\mathbb{C}^{N\times M}$. Find the sparse solution $\Vec{x}_{\ell_0}\in\mathbb{C}^M$ for given sparsity order $K_0\in\mathbb{N}$ such that the squared data residuals are minimal,
\be
\Vec{x}_{\ell_0}=\argmin_{\Vec{x}} \left\| \Vec{y} - \Mat{A} \Vec{x}  \right\|_2^2 \quad
    \text{subject to} \quad \| \Vec{x} \|_0 \leq K_0 ~,
\label{eq:l0_prob} \tag{P0}
\ee
where $\|\cdot\|_p$ denotes the $\ell_p$-norm.
The problem (\ref{eq:l0_prob}) is known as $\ell_0$-reconstruction. It is \emph{non-convex} and hard to solve \cite{Needell-Tropp-2009}. 
Therefore, the $\ell_0$-constraint in (\ref{eq:l0_prob}) is commonly relaxed to an $\ell_1$  constraint which renders the problem (\ref{eq:l1_prob}) to be \emph{convex}. 
Further, a matrix $\Mat{D}$ is introduced in the formulation of the constraint which gives flexibility in the problem definition. 
Let the number of rows of $\Mat{D}$ be arbitrary at first. 
In Sec. \ref{sec:dualproblem} suitable restrictions on $\Mat{D}$ are imposed where needed.
Several variants are discussed in \cite{Tibshirani2011}.
This gives
\be
\Vec{x}_{\ell_1}=\argmin_{\Vec{x}} \left\| \Vec{y} - \Mat{A} \Vec{x}  \right\|_2^2  \quad  \text{subject to} \quad
 \|\Mat{D} \Vec{x} \|_1 \leq \varepsilon~.
\label{eq:l1_prob} \tag{P1}
\ee
In the following, problem (\ref{eq:l1_prob}) is referred to as the complex-valued generalized LASSO problem.
 Incorporating the $\ell_1$  norm constraint into the objective function results in the equivalent formulation (\ref{eq:lasso_prob}),
\be 
\Vec{x}_{\ell_1}=\argmin\limits_{\Vec{x}} \left( \| \Vec{y} - \Mat{A} \Vec{x} \|_2^2 + \mu \| \Mat{D}\Vec{x} \|_1 \right)~. \tag{P1\mbox{$'$}}
\label{eq:lasso_prob}
\ee
The equivalence of (P0) and (P1') requires suitable conditions to be satisfied such as the restricted isometry property (RIP) condition or mutual coherence condition imposed on $\Mat{A}$, cf. \cite{DonohoElad03,CandesTao05,candes06}. Under such condition,  the problems (\ref{eq:l0_prob}) and (\ref{eq:lasso_prob}) yield the same sparsity order, $K_0=K$ with $K=\|\Vec{x}_{\ell_1}\|_0$, if the regularization parameter $\mu$ in (\ref{eq:lasso_prob}) is suitably chosen. The algorithms of Section \ref{se:algo} 
calculate suitable regularization parameters in this sense.

\section{Dual problem to the generalized LASSO}
\label{sec:dualproblem}

The generalized LASSO problem \cite{Tibshirani2011} is written in constraint form, all vectors and matrices are assumed to be complex-valued.
The following discussion is valid for arbitrary $N,M\in\mathbb{N}$: both the over-determined and the under-determined cases are included.
{
Following \cite{boyd,Kim-Boyd-Gorinevsky2009}, a  vector
 $\Vec{z}\in\mathbb{C}^M$ and
an  equality constraint $\Vec{z} = \Mat{D}\Vec{x}$ are introduced to obtain the equivalent problem
\be
  \label{eq:generalized-lasso-rewritten}
  \min\limits_{\Vec{x},\Vec{z}} \left(
  \| \Vec{y} - \Mat{A} \Vec{x} \|_2^2 + \mu \| \Vec{z} \|_1 \right)
 \quad  \text{subject to }\quad \Vec{z} = \Mat{D}\Vec{x}  ~.
\ee
The complex-valued dual vector $\Vec{u}=(u_1,\ldots,u_M)^T$  is introduced and associated with this equality constraint. The corresponding Lagrangian is }
\begin{eqnarray} \label{eq:tibshirani-lagrangian}
  \mathcal{L}(\Vec{x},\Vec{z},\Vec{u}) \!\!\!\!\!\!
  &=& \!\!  \!\!  \| \Vec{y} - \Mat{A} \Vec{x} \|_2^2 + \mu \| \Vec{z} \|_1+ 
     \Re\left[\Vec{u}^H(\Mat{D}\Vec{x}-\Vec{z})\right] 
     \\
&=&  \mathcal{L}_1(\Vec{x},\Vec{u})  +  \mathcal{L}_2(\Vec{z},\Vec{u}).
\end{eqnarray}
To derive the dual problem, the Lagrangian is minimized over $\Vec{x}$ and $\Vec{z}$. 
The terms involving $\Vec{x}$ are 
\be\label{eq:terms-in-xk}
 \mathcal{L}_1(\Vec{x},\Vec{u}) = \|\Vec{y}-\Mat{A}\Vec{x}\|^2_2 + \Re\left(\Vec{u}^H\Mat{D}\Vec{x} \right).
\ee
{The terms in (\ref{eq:tibshirani-lagrangian}) involving $\Vec{z}$ are
\be\label{eq:L2-involving-z}
  \mathcal{L}_2(\Vec{z},\Vec{u}) =  \mu \| \Vec{z} \|_1 - \Re(\Vec{u}^H\Vec{z})~.
\ee
}
The value $\hat{\Vec{x}}$ minimizing (\ref{eq:terms-in-xk}) is found by differentiation,
$\partial \mathcal{L}_1/\partial \Vec{x}=0$. 
This gives
\be\label{eq:Tibshirani-Eq34-Erich}
   \Mat{D}^H\Vec{u} =  2\Mat{A}^H \left( \Vec{y}- \Mat{A} \hat{\Vec{x}} \right) 
\ee   
whereby
\be\label{eq:Tibshirani-Eq34-Erich-followup}
   \Mat{A}^H \Mat{A} \hat{\Vec{x}} =  \Mat{A}^H \Vec{y} - \frac12 \Mat{D}^H\Vec{u}  ~.
\ee
If $ \Mat{D}^H\Vec{u}\in\mathop{\mathrm{span}}(\Mat{A}^H)$ the solution to
(\ref{eq:Tibshirani-Eq34-Erich-followup}) becomes,
\be\label{eq:Tibshirani-Eq34}
\hat{\Vec{x}} =  \underbrace{\Mat{A}^+ \Vec{y} +{\Vec{\xi}}}_{\hat{\Vec{x}}_{\mathrm{LS}}}-\frac12(\Mat{A}^H\Mat{A})^+ \Mat{D}^H\Vec{u}~,
\ee
where $(\cdot)^+$ denotes the Moore--Penrose pseudoinverse.
The Moore--Penrose pseudoinverse $\Mat{X}^+$ is defined and unique for all 
matrices $\Mat{X}$. 
In the following, we assume that $\Mat{A}$ has full row-rank and 
$\Mat{A}^+ = \Mat{A}^H (\Mat{A}\Mat{A}^H)^{-1}$ is a right-inverse \cite{GolubVanLoan}.
Here, $\Vec{\xi}\in\mathcal{N}(\Mat{A})$ is a nullspace term which enables $\hat{\Vec{x}}$ to deviate from the least norm solution $\Mat{A}^+\Vec{y}$. The nullspace $\mathcal{N}(\Mat{A})$ is $\{ \Vec{\xi}\in\mathbb{C}^M|\Mat{A}\Vec{\xi}=\Vec{0}\}$.
By identifying $\Vec{\xi}= \Vec{x}_{\ell_1}^{\mathrm{null}} $, we specialize (\ref{eq:Tibshirani-Eq34}) to the solution of (\ref{eq:lasso_prob}), 
\be\label{eq:Tibshirani-Eq34-specialized}
\Vec{x}_{\ell_1} = \Mat{A}^+ \Vec{y} + \Vec{x}_{\ell_1}^{\mathrm{null}} - \frac{_1}{^2} (\Mat{A}^H \Mat{A})^+ \Mat{D}^H \Vec{u}.
\ee
\begin{figure} 
	\centering
	\includegraphics[width=0.9\columnwidth]{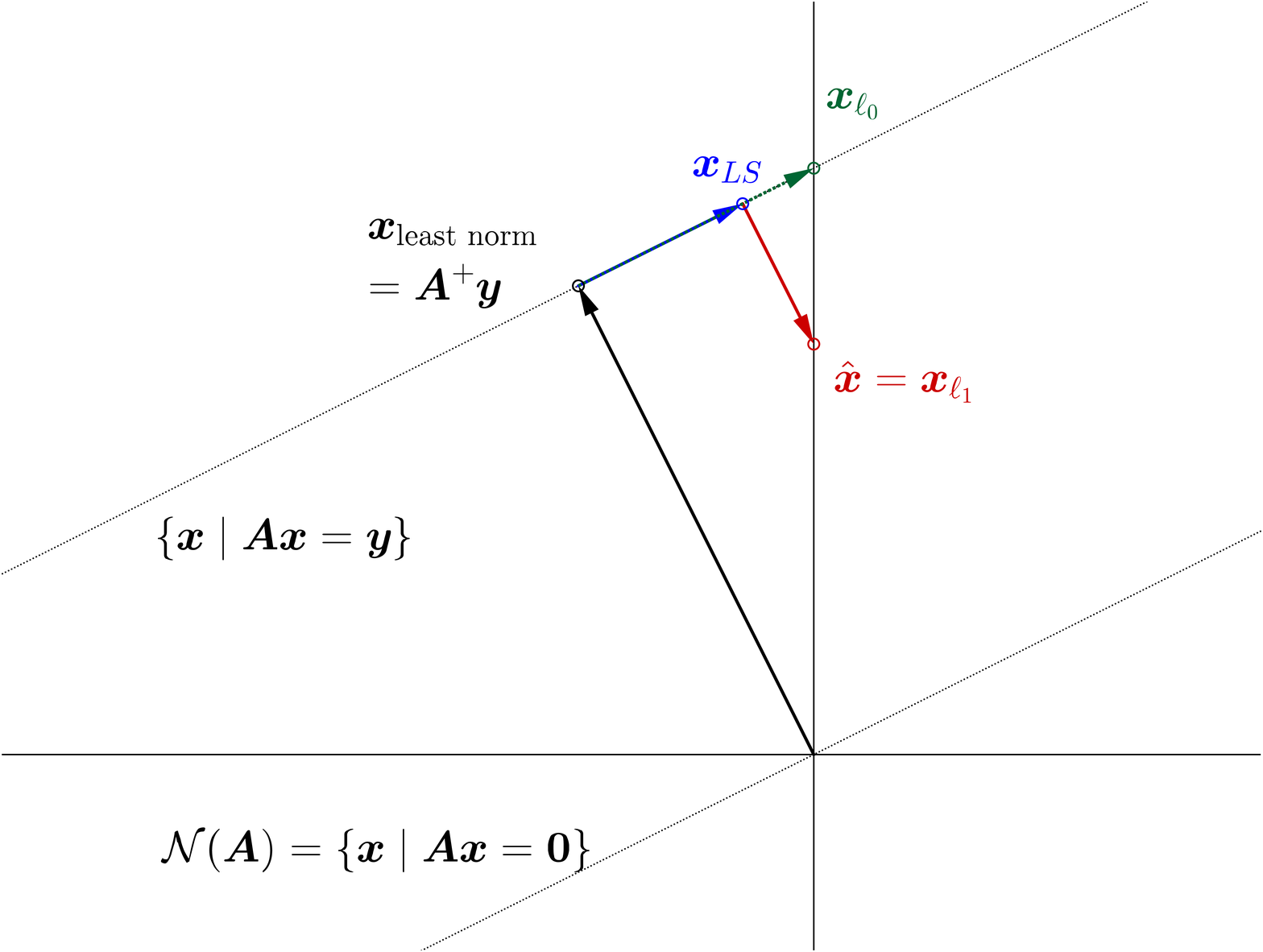}		
\caption{Sketch of the relations between the primal solution and the terms in (\ref{eq:Tibshirani-Eq34-specialized}): least norm solution $\Vec{x}_{\rm least\; norm}$, least squares solution $\Vec{x}_{\rm LS}$, and the sparse solutions $\Vec{x}_{\ell_0}$ , $\Vec{x}_{\ell_1}$. The nullspace term $\Vec{x}_{\ell_1}^{\mathrm{null}}$ is any vector along the line perpendicular to $\mathop{\mathrm{span}}(\Mat{A}^+)$. The red arrow represents the last term in (\ref{eq:Tibshirani-Eq34-specialized}) which is perpendicular to $\Vec{x}_{\ell_1}^{\mathrm{null}}$.}
	\label{fig:solution}
\end{figure}
\begin{figure} 
	\centering
      \includegraphics[width=0.9\columnwidth]{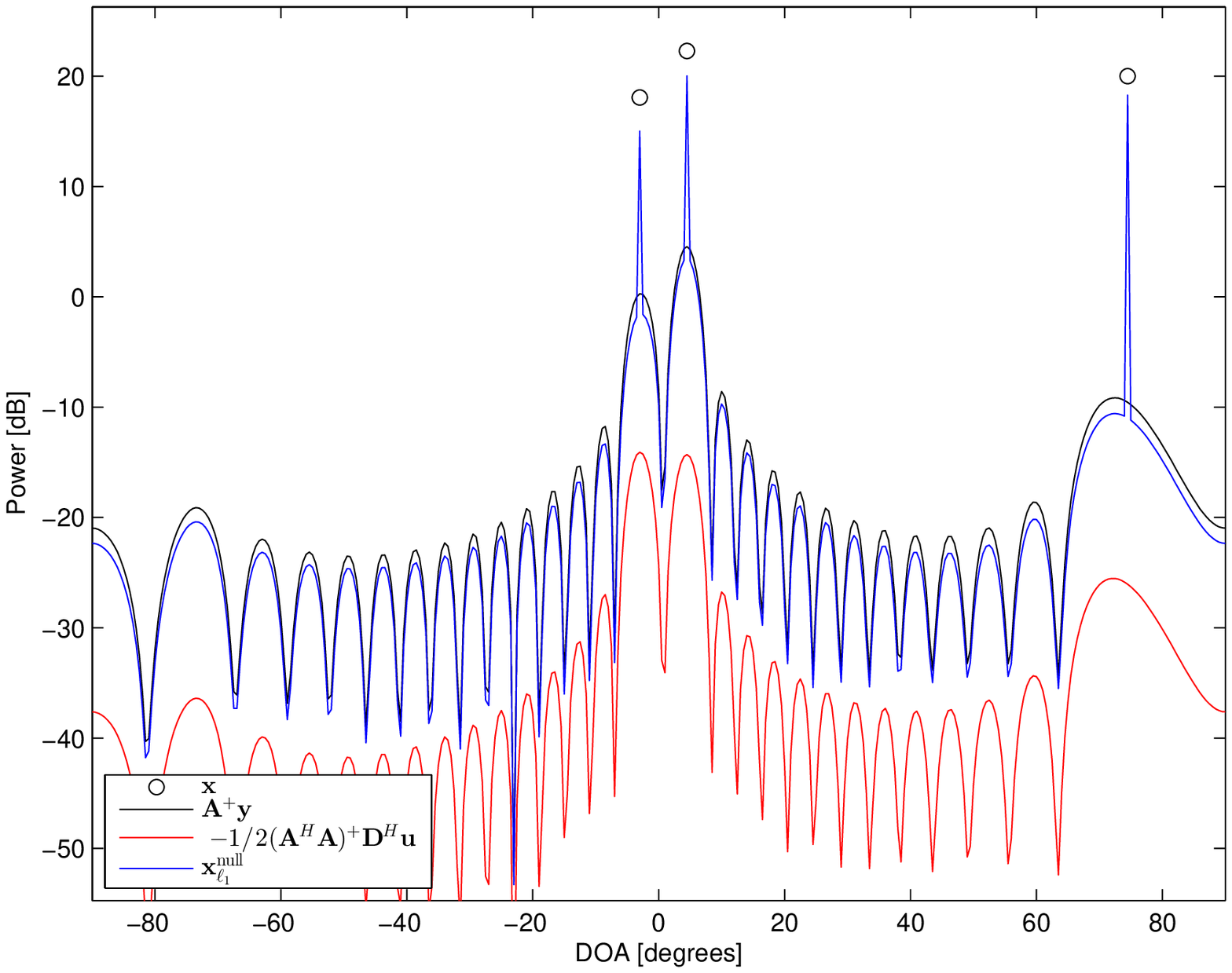}				
\caption{Numerical example solution terms in Eq.~(\ref{eq:Tibshirani-Eq34-specialized}) versus direction of arrival (DOA).}
	\label{fig:nullspace}
\end{figure}

Thus, the solution to the generalized LASSO problem (\ref{eq:Tibshirani-Eq34-specialized}) consists of three terms, as
 illustrated  in Fig. \ref{fig:solution}. The first two terms are the least norm solution $\Mat{A}^+ \Vec{y}$ and the nullspace solution $\Vec{\xi}$ which together form the unconstrained least squares (LS) solution $\hat{\Vec{x}}_{\mathrm{LS}}$.
The third term in (\ref{eq:Tibshirani-Eq34-specialized}) is associated with the dual solution.
Fig. \ref{fig:nullspace} shows the three terms of (\ref{eq:Tibshirani-Eq34-specialized}) individually for a simple  array-processing scenario. The continuous angle $\theta$ is discretized uniformly in $[-90, 90]^{\circ}$ using 361 samples and the wavefield is observed by 30 sensors resulting in a complex-valued $30 \times 361$  $\Mat{A}$ matrix (see section \ref{sec:array_model}).
At those primal coordinates $m$ which correspond to directions of arrival at $-3^\circ$, $4.5^\circ$ and $74.5^\circ$ in Fig. \ref{fig:nullspace}, the three terms in (\ref{eq:Tibshirani-Eq34-specialized}) sum constructively giving a non-zero $x_m$ (``the $m$th source position is active''), while for  all other entries they  interfere  destructively. Constructive interference is illustrated in Fig. \ref{fig:solution} which is in constrast to the destructive interference when the three terms in (\ref{eq:Tibshirani-Eq34-specialized}) sum to zero. This is formulated rigorously in Corollary \ref{cor:boundary}.

We evaluate (\ref{eq:terms-in-xk}) at the minimizing solution $\hat{\Vec{x}}$ and express the result solely by the dual  $\Vec{u}$.  Firstly, we expand
\be\label{eq:peter1}
  \| \Vec{y} - \Mat{A}\hat{\Vec{x}} \|^2_2 = \| \Vec{y} \|_2^2 + \| \Mat{A} \hat{\Vec{x}} \|_2^2 -
      2 \Re \{ \Vec{y}^H \Mat{A} \hat{\Vec{x}} \}
\ee
Secondly using (\ref{eq:Tibshirani-Eq34-Erich}), 
\begin{eqnarray}\label{eq:peter2}
  \Vec{u}^H \Mat{D} \hat{\Vec{x}} =  (\Mat{D}^H \Vec{u})^H \hat{\Vec{x}} &=& 2(\Vec{y} - \Mat{A} \hat{\Vec{x}})^H \Mat{A} \hat{\Vec{x}} \nonumber \\ 
   &=& 2\Vec{y}^H \Mat{A} \hat{\Vec{x}} - 2 \| \Mat{A} \hat{\Vec{x}} \|_2^2
    \label{eq:Erich11}
\end{eqnarray}
Adding Eq.(\ref{eq:peter1}) and the real part of (\ref{eq:Erich11}) gives
\begin{eqnarray}
      \mathcal{L}_1(\Vec{\hat x},\Vec{u}) &=&  \| \Vec{y} \|_2^2 - \| \Mat{A} \hat{\Vec{x}} \|_2^2 \nonumber\\
         &=& \Vec{y}^H \Vec{y} - \| \tilde{\Vec{y}} - \tilde{\Mat{D}}^H\Vec{u} \|_2^2 ~,
\end{eqnarray}
where we used (\ref{eq:Tibshirani-Eq34}) which assumes $\Mat{D}^H\Vec{u}\in\mathop{\mathrm{span}}(\Mat{A}^H)$ and introduced the abbreviations
\begin{eqnarray}
     \tilde{\Mat{D}} &=& \frac{_1}{^2} \Mat{D} \Mat{A}^+,  \label{eq:abbrev1} \\
     \tilde{\Vec{y}} &=& \Mat{P}_{\!_{\!\Mat{A}}} \,\Vec{y}~,\qquad
\text{with}~\Mat{P}_{\!_{\!\Mat{A}}}=\Mat{A} \Mat{A}^+ 
\end{eqnarray}
Due to the fundamental theorem of linear algebra, for an arbitrary vector $\Vec{v} \in \mathop{\mathrm{span}}(\Mat{A}^H)$ can be formulated as $\Mat{U}^H \Vec{v} = \Vec{0}$, where  
$\Mat{U}$ is a unitary basis of the null space $\mathcal{N}(\Mat{A})$. 
With $\Vec{v} = \Mat{D}^H \Vec{u}$, this becomes $(\Mat{DU})^H\Vec{u}=\Vec{0}$, resulting in
\be
 \inf\limits_{\Vec{x}} \mathcal{L}_1(\Vec{x},\Vec{u}) = \left\{ \begin{array}{ll}
 \Vec{y}^H\Vec{y} - \| \tilde{\Vec{y}} - \tilde{\Mat{D}}^H\Vec{u}\|_2^2~, & \hspace{-1ex}
 \text{if } (\Mat{D}\Mat{U})^H\Vec{u}=\Vec{0}, \\
   \hspace{-1ex} -\infty, &  \text{otherwise.}
  \end{array} \right. 
 \label{eq:almost-dual-function2}
\ee
Next (\ref{eq:L2-involving-z}) is minimized with respect to $\Vec{z}$, see Appendix A,
\begin{eqnarray}
    \inf\limits_{\Vec{z}}   \mathcal{L}_2(\Vec{z},\Vec{u})  &=& 
   \left\{ \begin{array}{ll}
    \hspace{-0.8ex} 0, &   \hspace{-0.8ex}\text{if } \|\Vec{u} \|_{\infty} \le \mu \\
    \hspace{-0.8ex} -\infty, &  \hspace{-0.8ex}\text{otherwise.} \end{array} \right.
    \label{eq:auxiliary-lemma}
\end{eqnarray}
Combining (\ref{eq:almost-dual-function2}) and (\ref{eq:auxiliary-lemma}),  the dual problem to the generalized LASSO (\ref{eq:l1_prob}) is,
\alpheqn
\begin{eqnarray} \label{eq:tibshirani-dual-problem35-min}
   \max\limits_{\Vec{u}\in\mathbb{C}^M} \Vec{y}^H\Vec{y} -\| \tilde{\Vec{y}} - \tilde{\Mat{D}}^H\Vec{u}\|_2^2  \\ \text{subject to }  \hspace{10ex}\|\Vec{u} \|_{\infty} \le \hspace{0.5ex}\mu, && \label{eq:box-constraint} \label{eq:boundary}\\
   (\Mat{D}\Mat{U})^H \Vec{u}\hspace{2.2ex} = \hspace{1ex}\Vec{0}. && \label{eq:row-space-constraint}
\end{eqnarray}
\reseteqn
Equation (\ref{eq:Tibshirani-Eq34-Erich}) is solvable for $\Vec{u}$ if the row space constraint (\ref{eq:row-space-constraint}) is fulfilled.
In this case, solving (\ref{eq:Tibshirani-Eq34-Erich}) directly gives 
\begin{theorem}\label{theo:dual}
{If $\Mat{D}$ is non-singular, }the dual vector $\Vec{u}$  is the output of a weighted matched filter acting on the vector of residuals, i.e.
\be\label{eq:2nd_time_mystery}
 \Vec{u} = 2 \Mat{D}^{-H} \Mat{A}^H ( \Vec{y} - \Mat{A} {\Vec{x}_{\ell_1}} )~,
\ee
where $\Vec{x}_{\ell_1}$ is the generalized LASSO solution (\ref{eq:lasso_prob}).
\end{theorem}
The dual vector $\Vec{u}$ gives an indication of the sensitivity of the primal solution to small changes in the constraints of the primal problem (cf.~\cite{boyd}: Sec.~5.6).
For the real-valued case the solution to (\ref{eq:lasso_prob}) is more easily constructed and better understood via the dual problem \cite{Tibshirani2011}.
Theorem \ref{theo:dual} asserts a linear one-to-one relation between the corresponding dual and primal solution vectors also in the complex-valued case. Thus, any results formulated in the primal domain are readily applicable in the dual domain.
This allows a more fundamental interpretation of sequential Bayesian approaches to density evolution for sparse source reconstruction \cite{Panahi-Viberg2011,MGPV2013}: they can be rewritten in a form that shows that they solve a generalized complex-valued LASSO problem and its dual. 
It turns out that the posterior probability density is strongly related to the dual solution 
\cite{MGPV2013,ZoechmanCoSeRa2015}. 

The following corollaries clarify useful element-wise relations between the primal and dual solutions: Corollary \ref{cor:boundary} relates the  \emph{magnitudes} of the corresponding primal and dual coordinates. Further, Corollary \ref{cor:phase} certifies what conditions on $\Mat{D}$ are sufficient for guaranteeing that the \emph{phase angles} of the corresponding primal and dual coordinates are equal.
Finally, Corollary \ref{cor:piecewise-linear} states that \emph{both} the primal and the dual solutions to (\ref{eq:lasso_prob}) are piecewise linear in the regularization parameter $\mu$.
\begin{cor} \label{cor:boundary}
For a diagonal matrix $\Mat{D}$ with real-valued positive diagonal entries, we conclude:
If the $m$th primal coordinate is active, i.e. $x_{\ell_1,m}\neq 0$ then the box constraint (\ref{eq:box-constraint}) is 
 tight in the $m$th dual coordinate. Formally, 
\be  
x_{\ell_1,m}\neq 0 \quad
     \Rightarrow\quad |u_m| = \mu, \qquad (m=1,\ldots,M).
\label{eq:implication-only}
\ee
\end{cor}
\noindent The proof is given in Appendix B.

\begin{figure} 
\centering
	\includegraphics[width=0.8\columnwidth]{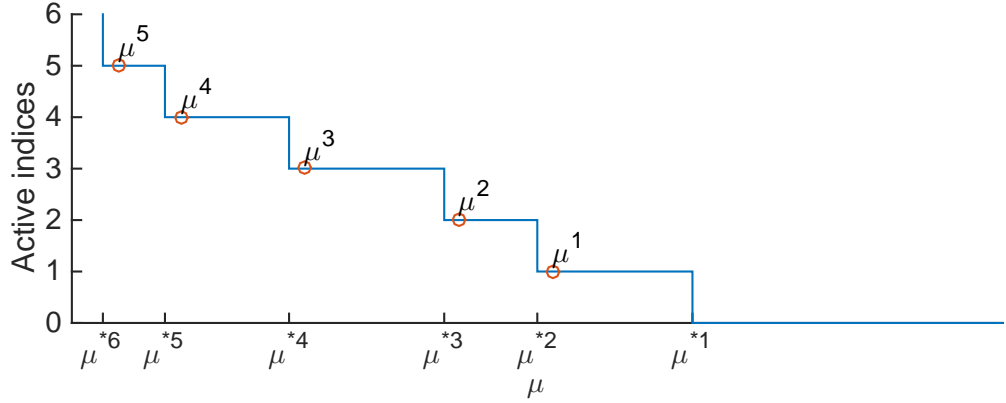}
\caption{Illustration of the LASSO path: Number of active indices versus the regularization parameter $\mu$. Increments in the active set occur at $\mu^{*p}$.} \label{fig:lassopath}
\end{figure}

Thus, the $m$th dual coordinate hits the boundary
as the $m$th primal coordinate becomes active. 
Conversely, when the bound on $|u_m|$ is loose (i.e. the constraint on $u_m$ is \emph{inactive}), the corresponding primal variable $x_m$ is zero (the $m$th primal coordinate is \emph{inactive}).
The active set $\mathcal{M}$ is 
\be \label{eq:active-set}
\mathcal{M} \;=\; \left\{m \, {\big|} \, x_{\ell_1,m}\neq 0 \right\} \; \subseteq \; \left\{m \, {\big|} \, |u_m| = \mu  \right\}
  \;=\; \mathcal{U}~.
\ee
Here, we have also defined the dual active set $\mathcal{U}$ which is a superset of $\mathcal{M}$ in general. This is due to Corollary \ref{cor:boundary} which states an implication in (\ref{eq:implication-only}) only, but not an equivalence.
The active set $\mathcal{M}$ implicitly depends on the choice of $\mu$ in problem (\ref{eq:lasso_prob}). 
Let $\mathcal{M}$ contain exactly $K$ indices,
\be
     \mathcal{M} = \{ m_1,\,m_2,\ldots,\,m_K\}.
\ee
The number of active indices versus $\mu$
 is illustrated in Fig.~\ref{fig:lassopath} \cite{Panahi-Viberg2012}.
Starting from a large choice of regularization parameter $\mu$ and then decreasing,
we observe incremental changes in the active set $\mathcal{M}$ at specific values $\mu^{*p}$ of the regularization parameter, i.e., the candidate points of the LASSO path \cite{Panahi-Viberg2012}.
The active set remains constant within the interval $\mu^{*p}>\mu>\mu^{*p+1}$.
By decreasing $\mu$, we enlarge the sets $\mathcal{M}$ and $\mathcal{U}$. By Eq.(\ref{eq:active-set}), we see that $\mathcal{U}$ may serve as a relaxation of the set of active indices $\mathcal{M}$. 
\begin{cor} \label{cor:phase}
If matrix $\Mat{D}$ is diagonal with real-valued positive diagonal entries, then the phase angles of the corresponding entries of the dual and primal solution vectors are equal.
\be 
\operatorname{arg}( {u}_m )= \operatorname{arg}( {x}_{\ell_1,m}), \quad \forall m \in \mathcal{M} \label{eq:phase_equal}
\ee
\end{cor}
\begin{cor} \label{cor:piecewise-linear}
The primal and the dual solutions to the complex-valued generalized LASSO problem (\ref{eq:lasso_prob}) are continuous and piecewise linear in the regularization parameter $\mu>0$. The changes in slope occur at those values for $\mu$ where the set of active indices $\mathcal{M}$ changes.
\end{cor}
\noindent The proofs for these corollaries are given in Appendix B.

\subsection{Relation to the $\ell_0$ solution}
It is now assumed that $\mathcal{M}$ defines the indices of the $K$ non-zero elements of the corresponding $\ell_0$ solution.
In other words: the $\ell_1$ and $\ell_0$ solutions share the same sparsity pattern. 
The $\ell_0 $ solution with sparsity order $K$ is then obtained by regressing the $K$ \emph{active} columns of $\Mat{A}$ to the data $\Vec{y}$ in the least-squares sense. Let 
\be \label{eq:col_match}
     \Mat{A}_{\mathcal{M}}  = [ \Vec{a}_{m_1},\,\Vec{a}_{m_2},\,\ldots,\Vec{a}_{m_K}]~,
\ee
where $\Vec{a}_m$ denotes the $m$th column of $\Mat{A}$.
The $\ell_0$ solution becomes (cf. Appendix C)
\be\label{eq:l0_sol}
      \Vec{x}_{\ell_0,\mathcal{M}} = \Mat{A}_\mathcal{M}^+ \Vec{y}~.
\ee
Here, $\Mat{A}_{\mathcal{M}}^+ = (\Mat{A}_{\mathcal{M}}^H\Mat{A}_{\mathcal{M}})^{-1}\Mat{A}_{\mathcal{M}}^H$ is the left inverse of $\Mat{A}_{\mathcal{M}}$.
By subtracting (\ref{eq:Tibshirani-Eq34-specialized}) from (\ref{eq:l0_sol}) and restricting the equations to the contracted basis $\Mat{A}_{\mathcal{M}}$ yields
\begin{eqnarray} \label{eq:relation}
\Mat{A}_{\mathcal{M}}(    \Vec{x}_{\ell_0,\mathcal{M}} - \Vec{x}_{\ell_1,\mathcal{M}}) &=& 
\frac{_1}{^2}\Mat{A}_{\mathcal{M}}\left( \Mat{A}^H_{\mathcal{M}} \Mat{A}_{\mathcal{M}} \right)^+ \Mat{D}_\mathcal{M}^H \Vec{u}_\mathcal{M} \nonumber \\
 &=& \tilde{\Mat{D}}^H_{\mathcal{M}} \mu e^{j \Vec{\theta}}~. 
\end{eqnarray}
In the image of $\Mat{A}$, the $\ell_0$-reconstruction problem (\ref{eq:l0_prob}) and the generalized LASSO (\ref{eq:lasso_prob}) coincide if the LASSO problem is pre-informed (prior knowledge) by setting $\Mat{D}_{mm}, \; m \in \mathcal{M}$ to zero. The prior knowledge is obtainable by an iterative re-weighting process \cite{Candes2008} or by a sequential algorithm on stationary sources \cite{MGPV2013}.


\section{Direction of Arrival Estimation}

For the numerical examples, we model a uniform linear array (ULA), which is described with its steering vectors representing the incident wave for each array element.

\subsection{Array Data Model} \label{sec:array_model}

Let $\Vec{x}=(x_{1},\ldots,x_{M})^{\mathrm T}$ 
be a vector of complex-valued source amplitudes. 
We observe time-sampled waveforms on an array of $N$ sensors 
which are stacked in the vector $ \Vec{y} $.
The following linear model for the narrowband sensor array data $ \Vec{y}$ at frequency $\omega$ is assumed, 
\be
   \Vec{y} = \Mat{A}\Vec{x} + \Vec{n}~.  \label{eq:linear-model2}
 \ee
The $m$th column of the transfer matrix $ \Mat{A}$ is the array steering vector $\Vec{a}_m$ for hypothetical waves from direction of arrival (DOA) $\theta_{m}$.
To simplify the analysis all columns are normalized such that their $\ell_2$ norm is one.
The transfer matrix $\Mat{A}$ is constructed by sampling all possible DOAs, but only very are active. 
Therefore, the dimension of $\Mat{A}$ is $N\times M$ with $N\ll M$ and $\Vec{x}$ is sparse.
The linear model (\ref{eq:linear-model2}) is underdetermined.

The $nm$th element of $\Mat{A}$ is modeled by 
\be\label{eq:A-elements}
  A_{nm} = \frac1{\sqrt{N}}\exp\left[\mathrm{j}(n-1)\pi\sin\theta_m\right].
\ee 
Here $\theta_m = \frac{(m-1)180^{\circ}}{M}-90^{\circ}$ is the DOA of the $m$th hypothetical DOA to the $n$th array element.

The additive noise vector $\Vec{n}$ is assumed spatially uncorrelated and follows a zero-mean complex normal distribution with diagonal covariance matrix $\sigma^2\Mat{I}$.

Following a sparse signal reconstruction approach \cite{Tibshirani2011,MGPV2013}, this leads to minimizing the
generalized LASSO Lagrangian 
\be\label{eq:generalized-lasso0}
 \| \Vec{y} - \Mat{A} \Vec{x} \|_2^2 + 
           \mu \left\| \Mat{D}\Vec{x} \right\|_1 ~,
\ee
where the weighting matrix $\Mat{D}$ gives flexibility in the formulation of the penalization term in (\ref{eq:generalized-lasso0}). 
Prior knowledge about the source vector leads to various forms of $\Mat{D}$.
This provides a Bayesian framework for sequential sparse signal trackers \cite{Panahi-Viberg2011,MGPV2013,ZoechmanCoSeRa2015}.
Specific choices of $\Mat{D}$ encourage both sparsity of the source vector and  sparsity of their successive differences which is a means to express that the source vector is locally constant versus DOA \cite{Tibshirani-Saunders-Rosset2005}.
The minimization of (\ref{eq:generalized-lasso0}) constitutes a convex optimization problem.
Minimizing the generalized LASSO Lagrangian (\ref{eq:generalized-lasso0}) with respect to $\Vec{x}$ for given $\mu$, 
gives a sparse source estimate $\Vec{x}_{\ell_1}$.
If $\mathop{\mathrm{rank}}(\Mat{A}) < M$, (\ref{eq:generalized-lasso0}) is no longer strictly convex and may not have a unique solution, cf. \cite{Tibshirani2011}.

Earlier approaches formulated this as a (ordinary) LASSO problem \cite{Tibshirani1996,Fuchs2005,malioutov05} which is equivalent to (\ref{eq:generalized-lasso0}) when specializing to $\Mat{D}=\Mat{I}$.

\subsection{Basis coherence} \label{sec:coherence}
The following examples feature different levels of basis coherence in order to examine the solution's behavior. 
As described in \cite{Xenaki2014}, the \textit{basis coherence} is a measure of correlation between two steering vectors and defined as the inner product between atoms, i.e. the columns of $\Mat{A}$. The maximum of these inner products is called mutual coherence and is customarily used for performance guarantees of recovery algorithms. To state the difference formally:
\be 
\operatorname{coh}\left(\Vec{a_i}, \Vec{a_j}\right) = \Vec{a_i}^H \Vec{a_j}
\label{eq:coh}
\ee
\be
\operatorname{mutual\ coh}(\Mat{A})=\left\|\Mat{A}^H \Mat{A} - \Mat{I} \right\|_\infty
\ee
The mutual coherence is bounded between $0$ and $1$

The following noiseless example in Figs. \ref{fig:3sources} and \ref{fig:4sources} demonstrates the dual solution and the WMF output for $N=30$ and $M=361$.
In Fig. \ref{fig:3sources}, the LASSO with $\mu=1$ is solved for a scenario with three sources at DOA $-3^\circ,4.5^\circ,84.5^\circ$ and all sources have same power level and are in-phase (see Fig. \ref{fig:3sources}b),
whereas in Fig. \ref{fig:4sources}, an additional fourth source at $8^\circ$ is included. 

\subsubsection{Low basis coherence}

\begin{figure} 
	\centering
      \begin{picture}(84,75)(4,0)
	\put(-15,0){\includegraphics[width=0.8\columnwidth]{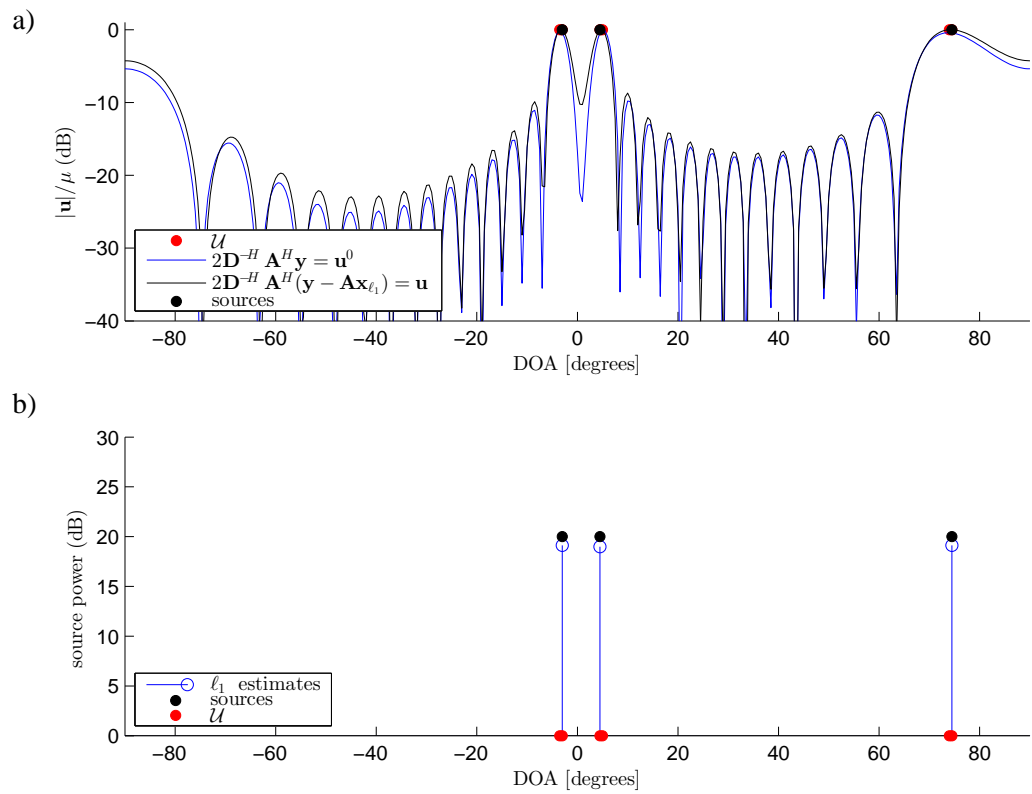}}
      \put(-20,102){a)}
      \put(-20,51){b)}
      \end{picture}			
      \caption{Dual (a) and primal (b) coordinates for 3 well separated sources with low basis coherence.}
	\label{fig:3sources}
\end{figure}

Figure \ref{fig:3sources} shows the performance when the steering vectors of the active sources have small basis coherence. The basis of source 1 is weakly coherent with source 2, $\operatorname{coh}\approx0.02$  using (\ref{eq:coh}).

Figure \ref{fig:3sources}a shows the normalized magnitude of the WMF (blue) and the normalized magnitude of the dual vector (black).
The dual active set $\mathcal{U}$ defined in (\ref{eq:active-set}) is depicted in 
red color.
 This figure shows that the true source parameters (DOA and power) are well estimated.
 It is also seen here that the behavior of the WMF closely resembles the magnitude of the dual vector and the WMF may be used as an approximation of the dual vector. This idea is further explored in Sec. \ref{se:algo}.

\unitlength1mm

\subsubsection{High basis coherence}

\begin{figure} 
	\centering
     \begin{picture}(84,78)(4,7)
	\put(-15,0){\includegraphics[width=0.8\columnwidth]{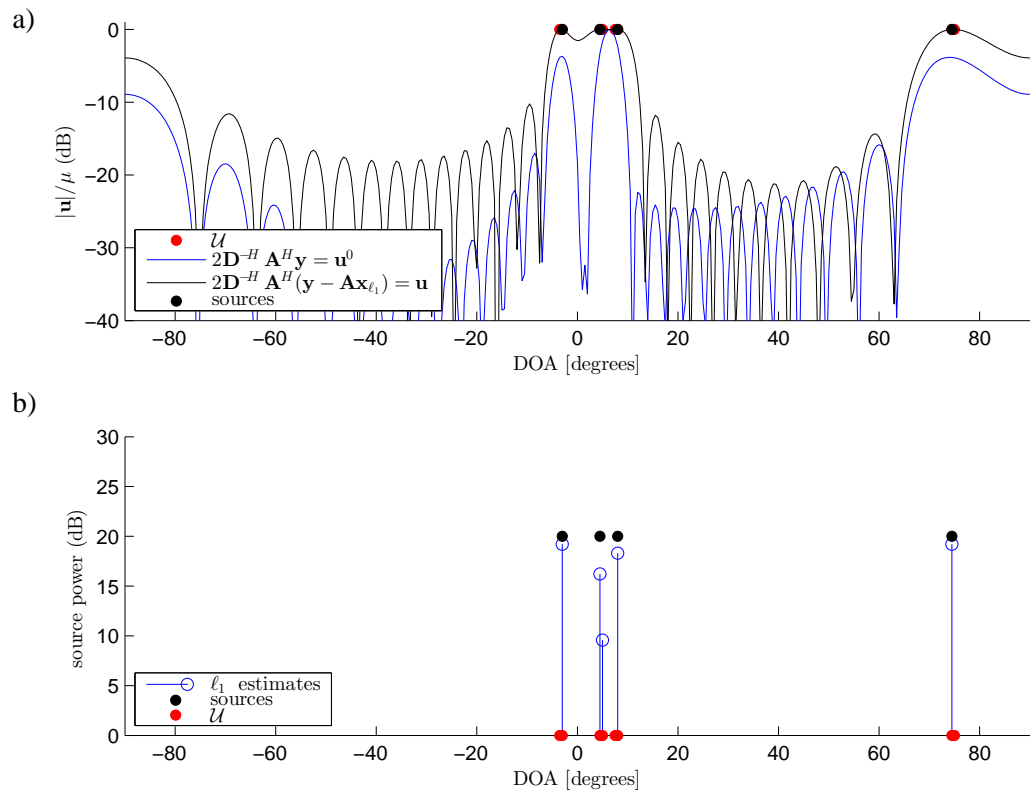}}
      \put(-20,102){a)}
      \put(-20,51){b)}
      \end{picture}			
      \caption{Dual (a) and primal (b) coordinates for 4 sources with higher basis coherence.}
\label{fig:4sources}
\end{figure}

Figure \ref{fig:4sources}a shows that the sources are not separable with the WMF, because the steering vectors belonging to source 2 and 3 are coherent, $\operatorname{coh}=0.61$  using (\ref{eq:coh}). The (generalized) LASSO approach is still capable of resolving all 4 sources. The DOA region defined by $\mathcal{U}$ is much broader around the nearby sources, allowing for spurious peaks close to the true DOA.
Figure \ref{fig:4sources}b shows that the true source locations (DOA) are still well estimated, but for the $2^{\rm nd}$ source from left, the power is split into two bins,
causing a poor source estimate.

\section{Solution path}

\begin{figure} 
      \centering
     \begin{picture}(84,65)(5,0)
	\put(-4,0){\includegraphics[width=0.8\columnwidth]{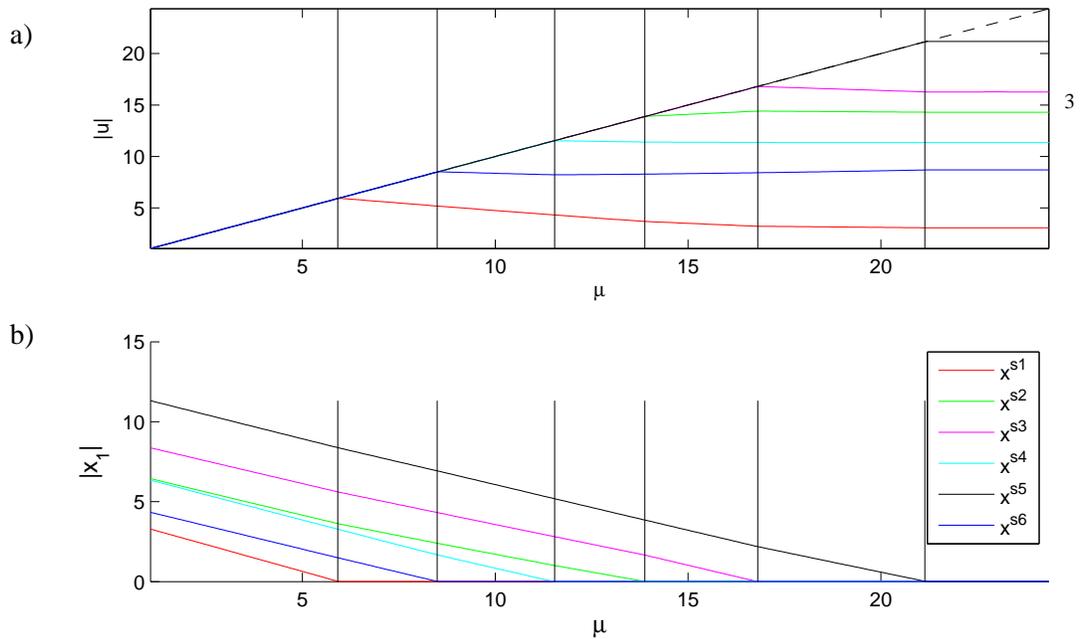}}
      \put(-12,80){a)}
      \put(-12,40){b)}
      \end{picture}			
\caption{Magnitudes of the solution paths versus $\mu$ for the simulation parameters in Table \ref{tab:src-params} and $\mathrm{SNR}=40\,$dB: (a) dual, and (b) primal vectors for the case of the \emph{complete basis}.} 
\label{fig:path6}
\end{figure}

\begin{figure} 
      \centering
  \begin{picture}(84,170)(5,0)
	\put(-6,0){\includegraphics[width=0.8\columnwidth]{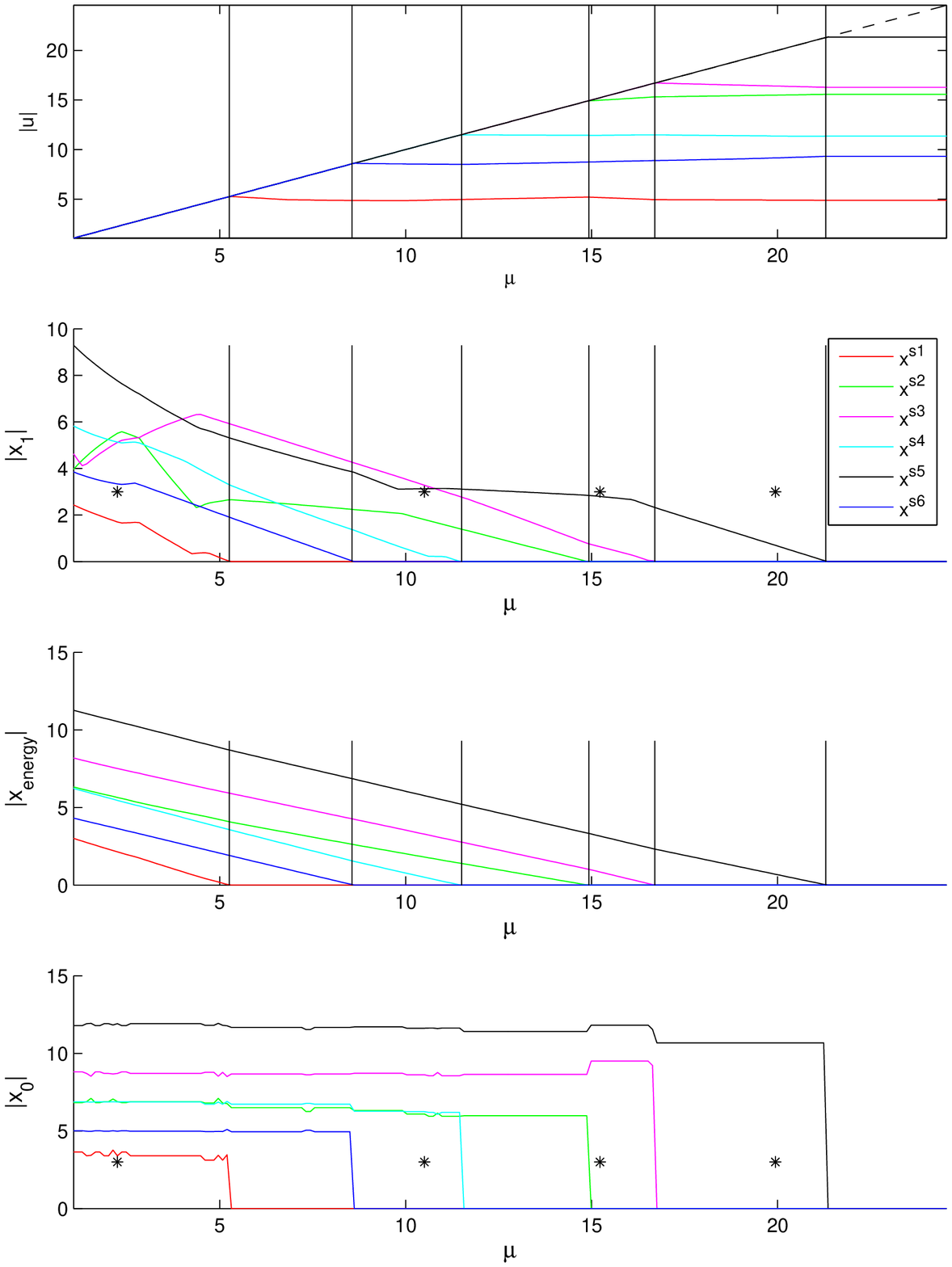}}
      \put(-12,170){a)}
      \put(-12,127){b)}
      \put(-12,85){c)}
      \put(-12,40){d)}
      \end{picture}	
\caption{Magnitudes of the solution paths versus $\mu$ for the simulation parameters in Table \ref{tab:src-params} and $\mathrm{SNR}=40\,$dB: (a) dual, and (b, c and d) primal vectors for the case of an 80-vector \emph{overcomplete basis}.
For the primal coordinates the peak within $\pm 2$ bins from the true bin is tracked based on (b) maximum (c) energy. The magnitudes of the corresponding elements of $\Vec{x}_{\ell_0}$ are shown in (d).} 
\label{fig:path5}
\end{figure}

\begin{figure} 
\centering
  \begin{picture}(84,150)(5,0)
	\put(-15,74){\includegraphics[width=0.8\columnwidth]{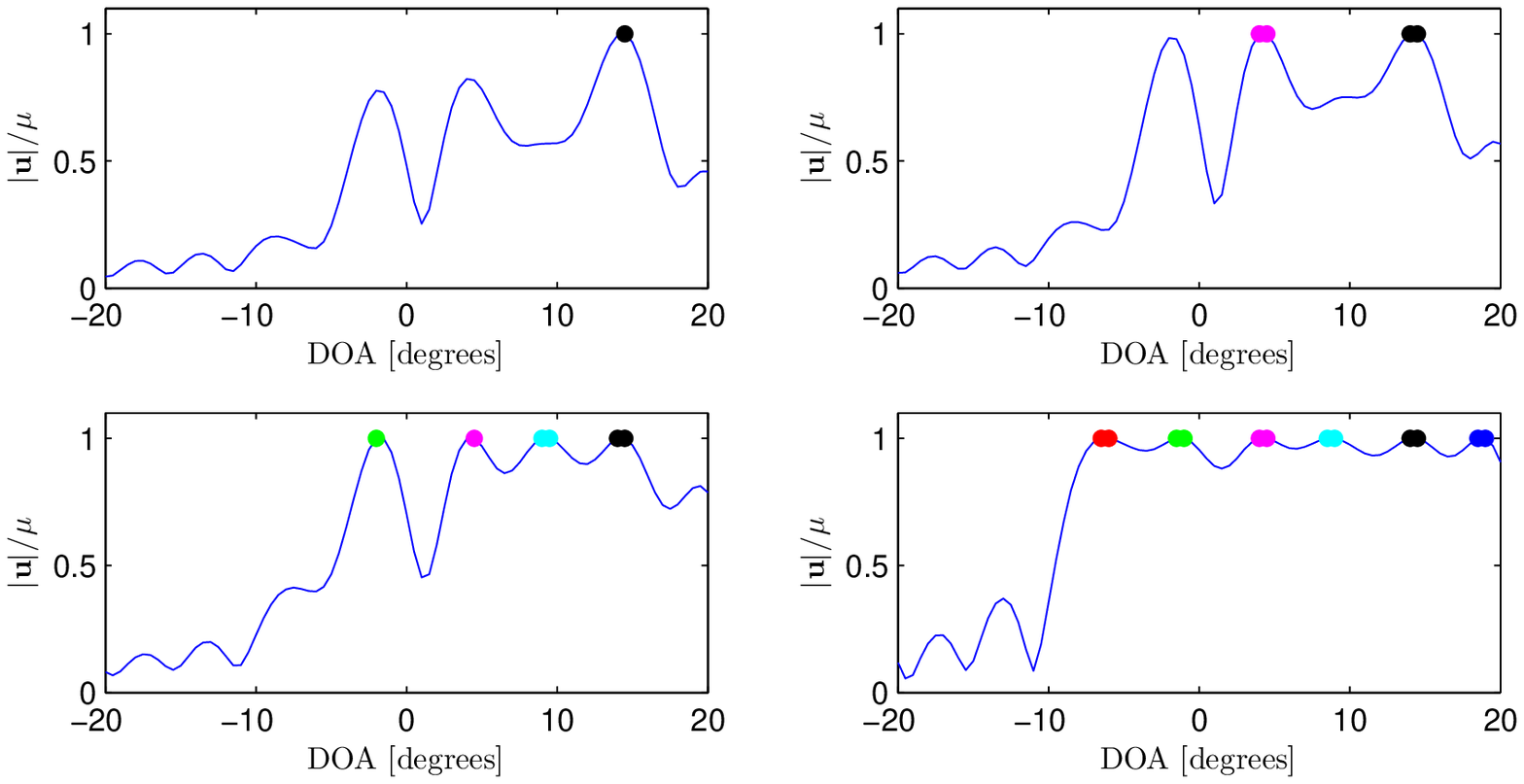}}
	\put(-15,0){\includegraphics[width=0.8\columnwidth]{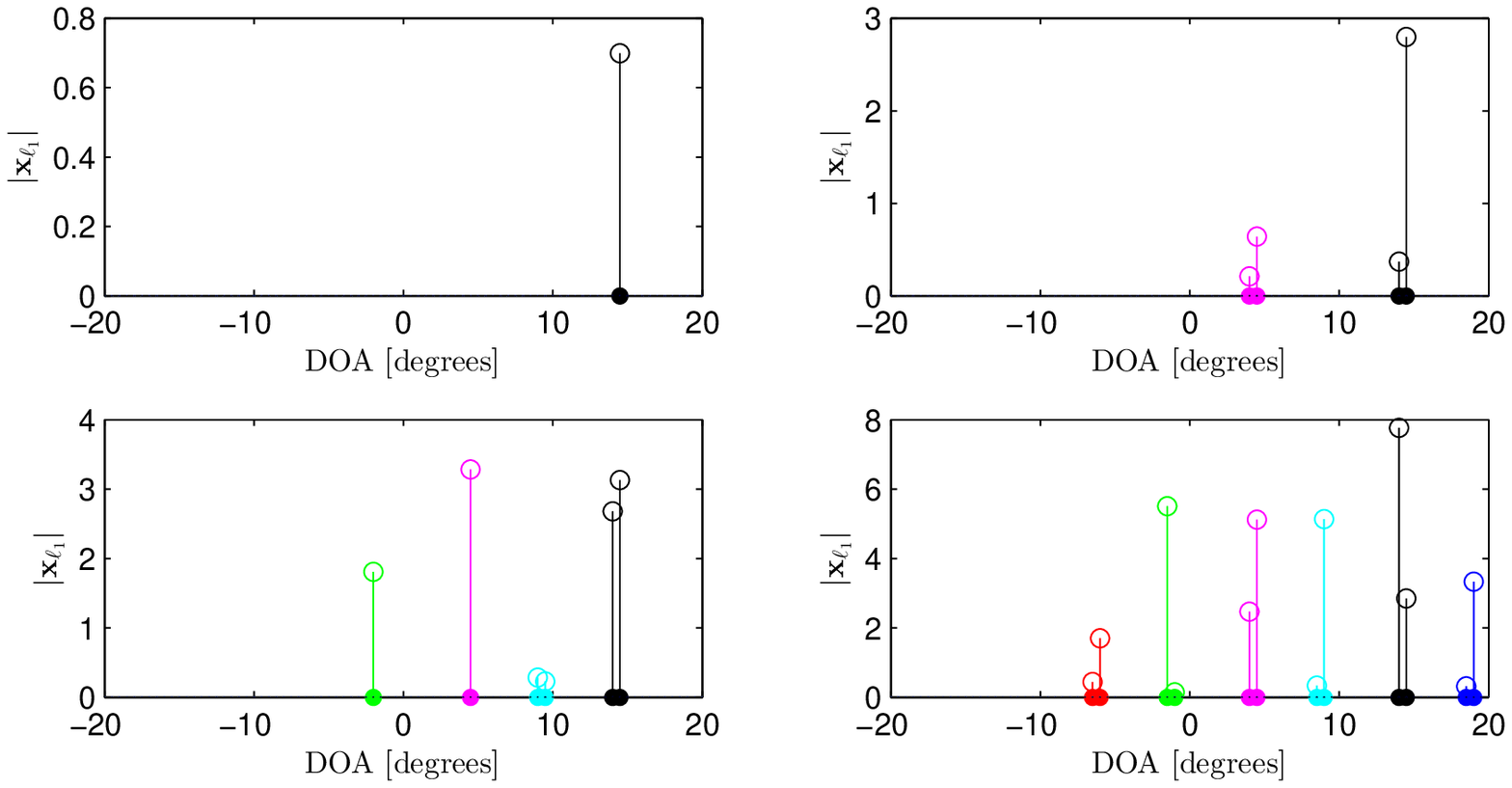}}
      \put(-20,140){\footnotesize a1)}        \put(51,140){\footnotesize b1)}
      \put(-20,104){\footnotesize c1)}        \put(51,104){\footnotesize d1)}
      \put(-20,65){\footnotesize a2)}         \put(51,65){\footnotesize b2)}
      \put(-20,30){\footnotesize c2)}         \put(51,30){\footnotesize d2)}
      \end{picture}	
\caption{Dual and primal coordinates at selected values of  $\mu $ for 81-vector overcomplete basis for $\mathrm{SNR}=40\,$dB.} \label{fig:path5dualprimal}
\end{figure}

\begin{figure} 
  \centering
  \begin{picture}(84,160)(5,0)
	\put(-15,0){\includegraphics[width=0.8\columnwidth]{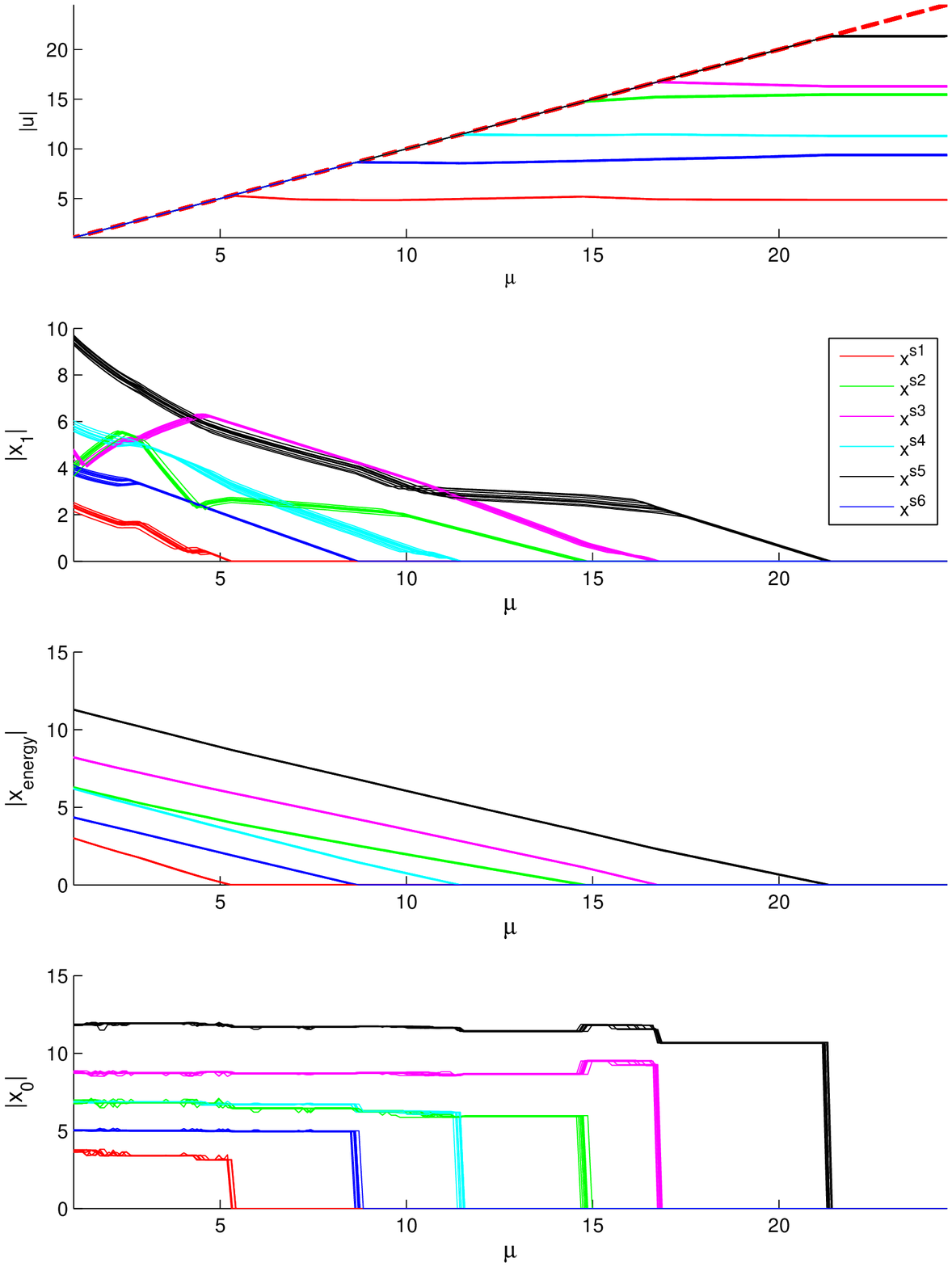}}
      \put(-20,167){a)}
      \put(-20,125){b)}
      \put(-20,83){c)}
      \put(-20,41){d)}
      \end{picture}	
\caption{For 10 noise realizations, magnitudes of the solution paths versus $\mu$ for the simulation parameters in Table \ref{tab:src-params} and $\mathrm{SNR}=40\,$dB: (a) dual, and (b, c and d) primal vectors for the case of an 80-vector \emph{overcomplete basis}.
For the primal coordinates the peak within $\pm 2$ bins from the true bin is tracked based on (b) maximum (c) energy. The magnitudes of the corresponding elements of $\Vec{x}_{\ell_0}$ are shown in (d).} 
 \label{fig:path5Noise40}
\end{figure}

\begin{figure} 
\centering
 \begin{picture}(84,160)(5,4)
	\put(-15,0){\includegraphics[width=0.8\columnwidth]{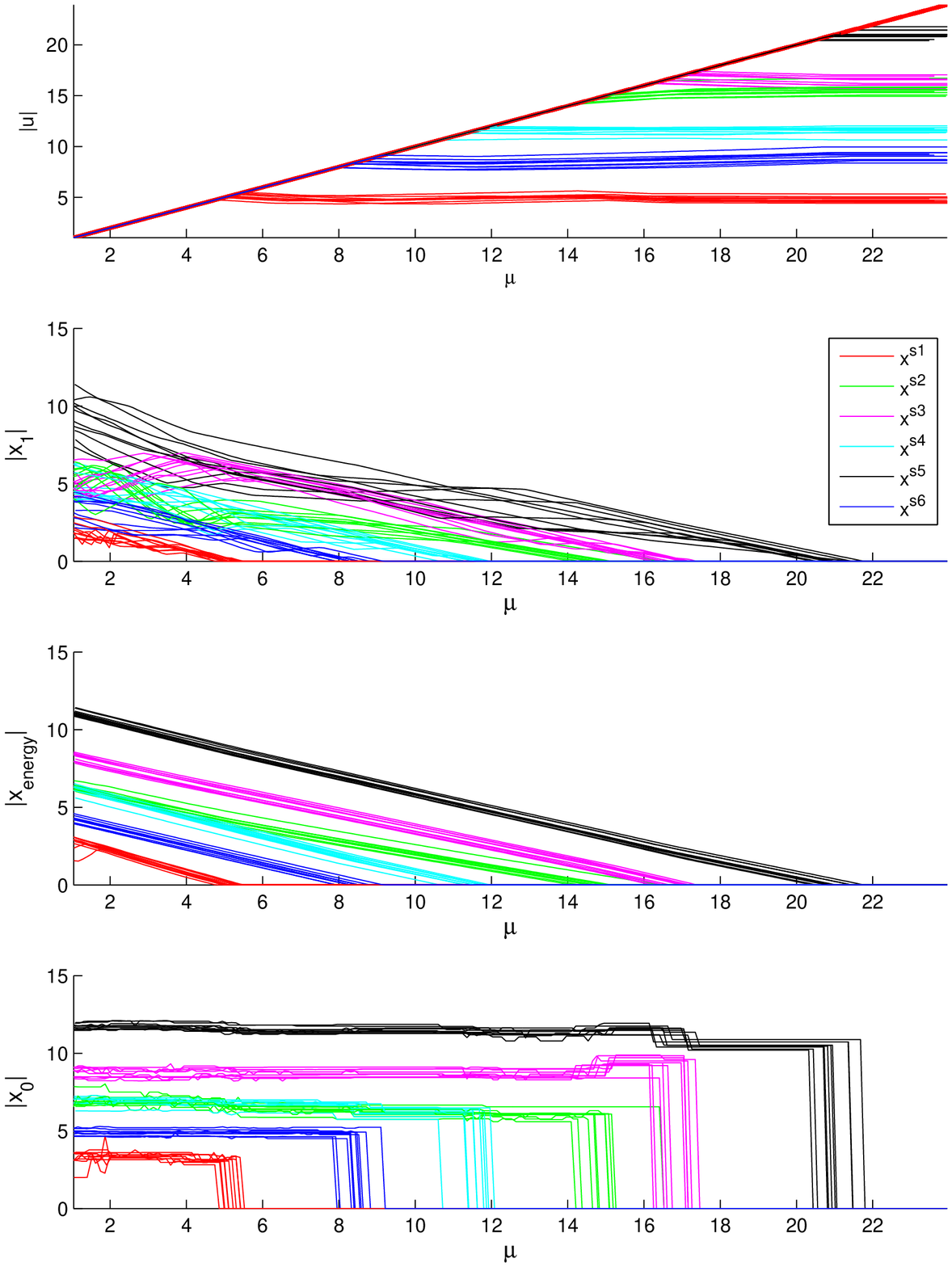}}
      \put(-20,167){a)}
      \put(-20,125){b)}
      \put(-20,83){c)}
      \put(-20,41){d)}
      \end{picture}	
\caption{As Fig. \ref{fig:path5Noise40}, but with $\mathrm{SNR}=20\,$dB: } \label{fig:path5Noise20}
\end{figure}

The LASSO solution path \cite{Tibshirani2011,Panahi-Viberg2012} gives the primal and dual solution vector versus the regularization parameter $\mu$. The primal and dual trajectories are piece-wise smooth and related according to Corollaries \ref{cor:boundary}--\ref{cor:piecewise-linear}. The following figures show results from individual LASSO runs by varying $\mu$.

The problems (\ref{eq:l1_prob}) is complex-valued and the corresponding solution paths behave differently from what is described in Ref. \cite{Tibshirani2011}. In the following figures, only the magnitudes of the active primal coordinates and the corresponding dual coordinates are illustrated.
Note that Corollary \ref{cor:phase} guarantees that the phases of the active primary solution elements and their duals are identical and independent from $\mu$. 

Based on the observed solution paths, we notice that the hitting times (when $|u_m|=\mu$) of the dual coordinates (at lower $\mu$) are well predictable from the solution at higher $\mu$.

For the following simulations and Figs. \ref{fig:path6}--\ref{fig:path5Noise20}, the signal to noise ratio (SNR) is defined as
\be
\mathrm{SNR} = 10 \log_{10} \left( \mathrm{E}\| \Mat{A} \Vec{x} \|_2^2 \, / \, \mathrm{E}\| \Vec{n} \|_2^2 \right)\,\mathrm{dB}.
\ee 
$\mathrm{SNR}=40\,\mathrm{dB}$ in Figs. \ref{fig:path6}--\ref{fig:path5Noise40}, whereas $\mathrm{SNR}=20\,\mathrm{dB}$ in Fig. \ref{fig:path5Noise20}.

\subsection{Complete Basis}
First (Fig.~\ref{fig:path6})  discusses the dual and primal solution for a complete basis with $M=6$,  sparsity order $K=6$, and $N=30$ sensors linearly spaced with half wavelength spacing. 
This simulation scenario is not sparse and all steering vectors $\Vec{a}_m$ for $1\le m\le M$ will eventually be used to reconstruct the data for small $\mu$. 
The source parameters that are used in the simulation scenario are given in Table \ref{tab:src-params}. 

\begin{table}
\begin{center}
\begin{tabular}{rrr} \hline
No. & DOA ($^{\circ}$) & Power (lin.) \\ \hline
1 & $-$6.0 & 4.0 \\
2 & $-$1.0 & 7.0 \\
3 & 4.0 & 9.0 \\
4 & 9.0 &  7.0 \\
5 & 14.0 & 12.0 \\
6 & 19.0 & 5.0 \\ \hline
\end{tabular}
\caption{Source parameters for simulation scenario}
\label{tab:src-params}
\end{center}
\end{table}

We discuss the solution paths in Figs. \ref{fig:path6}--\ref{fig:path5Noise20} from right ($\mu = \infty$) to left ($\mu = 0$). Initially all dual solution paths are horizontal (slope = 0), since the primal solution  $\Vec{x}_{\ell_1}=\Vec{0}$ 
for $\mu > 2 \| \Mat{D}^{-H} \Mat{A}^H\Vec{y}\|_{\infty}$.
In this strongly penalized regime, the dual vector is the output of 
the WMF $\Vec{u}= 2\Mat{D}^{-H} \Mat{A}^H\Vec{y}$ which does not depend on $\mu$. 

At the point $\mu^1 = 2 \| \Mat{D}^{-H} \Mat{A}^H\Vec{y}\|_{\infty}$ the first dual coordinate hits the boundary (\ref{eq:boundary}). This occurs at $\mu^1=21$ in Fig.~\ref{fig:path6}a and the corresponding primal coordinate becomes active. As long the active set $\mathcal{M}$ does not change, the magnitude of the corresponding dual coordinate is $\mu$, due to Corollary \ref{cor:boundary}. The remaining dual coordinates change slope relative to the basis coherence level of the active set. 

As $\mu $ decreases, the source magnitudes at the primal active indices increase since the 
$\ell_1$-constraint in (\ref{eq:lasso_prob}) becomes less important, see Fig.~\ref{fig:path6}b.
The second source will become active when the next dual coordinate hits the boundary (at $\mu^1=17$ in Fig.~\ref{fig:path6}). 

When the active set is constant, the primary and dual solution is piecewise linear with $\mu$, as proved in Corollary \ref{cor:piecewise-linear}. The changes in slope are quite gentle, as shown for the example in Fig. \ref{fig:path6} .
Finally, at $\mu=0$ the problem (\ref{eq:lasso_prob}) degenerates to an 
unconstrained (underdetermined) least squares problem.  
Its primal solution $\hat{\Vec{x}}= \hat{\Vec{x}}_{\mathrm{LS}}$, see (\ref{eq:Tibshirani-Eq34}), is not unique and the dual vector is trivial, $\Vec{u}= \Vec{0}$.

\subsection{Overcomplete Basis}

We now enlarge the basis to $M=81$ with hypothetical source locations $\theta_m \in [-20^\circ, \, 20^\circ]$ with $0.5^\circ$ spacing, and all other parameters as before.
The solution is now sparse.

The LASSO path \cite{Panahi-Viberg2012} is illustrated in  Fig.~\ref{fig:path5} where we expect the 
source location estimate within $\pm 2 $ bins from the true source location.  
The dual  Fig.~\ref{fig:path5}a appears to be quite similar to  Fig.~\ref{fig:path6}a.

Corollary \ref{cor:piecewise-linear} gives that the primary solution should change linearly, as demonstrated for  the complete basis in Fig.~\ref{fig:path6}b. Here we explain why this is not the case for  the overcomplete basis
 primary solution in Fig.~\ref{fig:path5}b . This is understood by examining the full solution at selected values of $\mu$ (asterisk (*) in  Fig.~\ref{fig:path5}).
At $\mu=20$ just one solution is active, only the black source (source 5)  is active though one bin to the left, as shown in Fig.~\ref{fig:path5dualprimal}a2. The dual vector in Fig.~\ref{fig:path5dualprimal}a1--{fig:path5dualprimal}d1, has a broad maximum, explaining the sensitivity to offsets around the true DOA. The shape of this maximum is imposed by the dictionary; the more coherent the dictionary, the broader the maximum.
Between $\mu=16$ and $\mu=11$, the black source appears constant, this is because at large values the source is initially located in a neighboring  bin. As $\mu$ decreases, the correct bin receives more power, see Fig.~\ref{fig:path5dualprimal}b2 and Fig.~\ref{fig:path5dualprimal}c2 for $ \mu=15$ and $ \mu=10$, respectively. When it is stronger than the neighboring bin at $\mu\le11$, see Fig.~\ref{fig:path5dualprimal}d2, this source power starts increasing again. This trading in source power causes the fluctuations in  Fig.~\ref{fig:path5}b.

One way to correct for this fluctuation is to sum the coherent energy for all bins near a source, i.e., multiplying the source vector with the corresponding neighbor columns of $\Mat{A}$, which also touch the boundary (marked region in Fig. \ref{fig:path5dualprimal}) and then  compute the  energy based on the average received power at each sensor. This gives a steady rise in power as shown in   Fig.~\ref{fig:path5}c.

We motivate solving (\ref{eq:lasso_prob}) as a substitute for $\ell_0$-reconstruction (\ref{eq:l0_prob})---finding the active indexes of the $\ell_1$ solution, see  Fig.~\ref{fig:path5}d.
The $\ell_0$ primal can be found with the restricted basis and the value of the $\ell_1$ primal from (\ref{eq:Tibshirani-Eq34}), which depends on $\mu$, or by just solving (\ref{eq:l0_sol}).

To investigate the sensitivity to noise, 10 LASSO  paths are simulated for 10 noise realizations for  both  $\mathrm{SNR}=40\,$dB (Fig.~\ref{fig:path5Noise40}) and  $\mathrm{SNR}=20\,$dB (Fig.~\ref{fig:path5Noise20}).
The dual (Fig.~\ref{fig:path5Noise40}a and  Fig.~\ref{fig:path5Noise20}a), appears quite stable to noise, but the primal $|{\bf x}_{\ell 1}|$ (Figs. \ref{fig:path5Noise40}b and   \ref{fig:path5Noise20}b) show quite large variation with noise. This is because the noise causes the active indexes to shift and thus the magnitude to vary. 
The mapping to energy   $|{\bf x}_{\rm energy}|$ (Figs. \ref{fig:path5Noise40}c and   \ref{fig:path5Noise20}c) or the $|{\bf x}_{\ell 0}|$ solution (Figs. \ref{fig:path5Noise40}d and   \ref{fig:path5Noise20}d)  makes the solution much more stable.

\section{Solution Algorithms}

\label{se:algo}

Motivated by Theorem \ref{theo:dual} and Corollary \ref{cor:boundary}, 
we propose the order-recursive algorithm in Table  \ref{algo:order-recursive_new} for approximately solving problem (\ref{eq:l0_prob}) by selecting a suitable regularization parameter $\mu$ in problem (\ref{eq:lasso_prob}), a faster iterative algorithm in Table \ref{algo:fast_new}, and a dual-based iterative algorithm in Table \ref{algo:dual-fast_new}.\\
\\
As shown by Theorem \ref{theo:dual}, the dual vector is evaluated by a WMF acting on the LASSO residuals. The components of the dual vector which hit the boundary, i.e. $|u_m|=\mu$, correspond to the active primal coordinates $|x_m|>0$. As $|u_m|=\mu$ constitutes a necessary condition, this condition is at least $|\mathcal{M}|$ times fulfilled. Informally, we  express this as: ``The dual vector must have $|\mathcal{M}|$ peaks of height $\mu$, where the shaping is defined by the dictionary $\Mat{A}$ and the weighting matrix $\Mat{D}$."

The key observation is the reverse relation. By knowing the peak magnitudes of the dual vector, one estimates the appropriate $\mu$-value to make $i$ peaks hit the boundary. We denote this regularization parameter value as $\mu^i$.  This is a necessary condition to obtain $i$ active sources. 

We define the $\mathop{\mathrm{peak}}(\Vec{u},i)$--function  which returns the $i^{\rm th}$ largest local peak in magnitude of the vector $\Vec{u}$.
A local peak  is defined as an element which is larger  than  its adjacent elements. The peak function can degenerate to a simple sorting function giving the $i$th largest value, this will cause slower convergence in the algorithms below. 
\begin{prop} \label{prop:peak}
Assuming all sources to be separated such that there is at least a single bin in between, the $\mathop{\mathrm{peak}}$ function relates the regularization parameter to the dual vector via
\be  \label{eq:fix_point}
\mu^i  ={\rm peak}\left(\big| \Vec u(\mu^i) \big|, i\right) ={\rm peak}\left( \Vec u(\mu^i), i\right) ~.
\ee
\end{prop}
Equation (\ref{eq:fix_point}) is a fixed-point equation for $\mu^i$ which is demanding to solve. Therefore we approximate (\ref{eq:fix_point}) with previously obtained dual vectors\footnote{For the first step, we define the WMF output as $\Vec{u}_0 = 2 \Mat{D}^{-H}\Mat{A}^H \Vec{y}$.}. 
At a potential new source position $n$, the dual vector is expanded as
\begin{eqnarray} 
u_n(\mu^i) &=& 
\frac{2}{{D}_{n, n}^*} \Vec{a}_{n}^H \bigg( \Vec{y}- \hspace{-1ex} \sum  \limits_{m \in \mathcal{M}_i} \Vec{a}_m \Vec{x}_{\ell_1,m}(\mu^i) \bigg)   \\ 
&\approx&  \frac{2} {{D}_{n, n}^*} \Vec{a}_{n}^H \bigg( \Vec{y}- \hspace{-2ex} \sum\limits_{m \in \mathcal{M}_{i-1}} \hspace{-2ex}\Vec{a}_m \Vec{x}_{\ell_1,m}(\mu^{i-1}) \bigg) \label{eq:approx_rec}  \\   
&\approx&  \frac{2} {{D}_{n, n}^*} \Vec{a}_{n}^H \bigg( \Vec{y}- \hspace{-2ex} \sum\limits_{m \in \mathcal{M}_{i-2}} \hspace{-2ex}\Vec{a}_m \Vec{x}_{\ell_1,m}(\mu^{i-2}) \bigg)  \\
&\vdots&  \nonumber \\   
&\approx&  \frac{2} {{D}_{n, n}^*} \Vec{a}_{n}^H  \Vec{y} \label{eq:approx_fast} 
\end{eqnarray}
The approximations used in (\ref{eq:approx_rec})--(\ref{eq:approx_fast}) are 
progressive. These approximations are good  if the steering vectors 
corresponding to the active set are sufficiently incoherent: 
$|\Vec{a}^H_n \Vec{a}_m| \approx 0$ for $n,m\in\mathcal{M}$. 
Eq. (\ref{eq:approx_fast}) corresponds to the conventional beamformer $\Mat{A}^H\Vec{y}$ for a single snapshot. In the solution algorithms, the approximations (\ref{eq:approx_rec})--(\ref{eq:approx_fast}) are used for the selection of the regularization parameter $\mu$ only, 
thus the peaks in the conventional beamformer do not correspond to 
 the $\Vec{x}_{\ell_1}$ solution.

Our simulations have shown that a significant speed-up achievable, so we named it fast-iterative algorithm, cf Section \ref{se:fast}. \\
\\
From the box constraint (\ref{eq:box-constraint}), the magnitude of the $i^{\rm th}$ peak in $\Vec{u}$ does not change much during the iteration over $i$: It is bounded by the difference in regularization parameter. For any $\mu^{i}<\mu^{i-1}$, we conclude from Corollary \ref{cor:boundary} and Proposition \ref{prop:peak}
\be
  \underbrace{\mathop{\mathrm{peak}}(\Vec{u}(\mu^{i-1}),i)}_{\le \mu^{i-1}} - \underbrace{\mathop{\mathrm{peak}}(\Vec{u}(\mu^{i}),i)}_{=\mu^{i}} \le \mu^{i-1} - \mu^{i}. \label{eq:prediction-error}
\ee
Thus, the magnitude of the $i^{\rm th}$ peak cannot change more than the corresponding change in the regularization parameter. The left hand side of (\ref{eq:prediction-error}) is interpretable as the prediction error of the regularization parameter and this shows that the prediction error is bounded.

Assuming our candidate point estimates ($\mu^{*1},\,\mu^{*2},\,\ldots$) are correct, we follow a path of regularization parameters $\mu^1,\,\mu^2,\,\ldots$ where $\mu^p$ is slightly higher than the lower end $\mu^{*p+1}$ of the regularization interval. Specifically, $\mu^p=(1-F)\mu^{*p}+F\mu^{*p+1}$ with $F<1$. For the numerical examples $F=0.9$ is used. 
This $F$ is chosen because the primal solution $\Vec{x}_{\ell_1}$ is closest to $\Vec{x}_{\ell_0}$ at the lower end of the interval.

In the following we focus on the order recursive algorithm, and indicate the differences to the other approaches.

\subsection{Recursive-In-Order algorithm}
\label{se:lassosol}

The recursive-in-order algorithm  in Table \ref{algo:order-recursive_new}  finds one source at a time as $\mu$ is lowered. To this purpose it employs an approximation of the height of the $i$th local peak given a solution with $(i-1)$ peaks. The underlying assumption is that the next source will become active at the location corresponding to the dual coordinate of the next peak. 
Equation (\ref{eq:approx_rec}) allows to approximate
\begin{eqnarray}
\mu^i  &=& {\rm peak}(\Vec{u}(\mu^{i}), i) \approx {\rm peak}(\Vec{u}(\mu^{i-1}), i)~. \label{eq:mu-approx}
\end{eqnarray}
This assumption is not universally valid as \delete{discussed in the last paragraph of this section.\\}
\delete{
The algorithm starts  with the all-zero solution $\Vec{x}_{\ell_1}^0 = \Vec{0} $ in line 1. This corresponds to a large value of 
$\mu=\mu^0>2\|\Mat{D}^{-H}\Mat{A}^H\Vec{y}\|_{\infty}$, but is not used directly. 
The first value $\mu^1$  is chosen based on the first peak in $\Vec{u}$ for $\Vec{x}_{\ell_1}=\Vec{0}$.
\delete{The $\ell_{\infty}$-norm is implemented by calling the $\mathop{\mathrm{peak}}(\Vec{u},p)$--function for $p=1$ in line 3.}

In lines 4--5 in Table \ref{algo:order-recursive}, the generalized LASSO problem (\ref{eq:lasso_prob}) is solved for $\mu=\mu^1$ and the corresponding  active set is detected by thresholding (at this point, the active set contains only a single active index). Then, by assuming that the active set is also valid for the solution to the 
$\ell_0$-problem (\ref{eq:l0_prob}), the solution $\Vec{x}_{\ell_0}^p$ is computed.
Next, we start an iteration loop with counter $i$ which reduces the regularization parameter $\mu$ 
by inserting the updated dual variable $\Vec{u}_1$ into line 3. 

The next choice of regularization parameter $\mu^2$ is selected in the interval $[\mu^{*2},\mu^{*3}[$ which we estimate by the $2^{\mathrm{nd}}$ and $3^{\mathrm{rd}}$ largest peaks of the dual variable. This is illustrated in Fig. \ref{fig:lassopath}.
This is implemented in line 3.
Then, we solve (\ref{eq:lasso_prob}) for $\mu = \mu^2<\mu^1$ and continue the iteration until the given sparsity order $K=K_0$ is reached. 

In line 5, the active set $\mathcal{M}$ is approximated by thresholding of the primal solution.

If the number of elements in $\mathcal{M}$ equals the loop counter $i$ then the basis is restricted to $\mathcal{M}$ and the corresponding $\ell_0$-solution is computed in lines 9--13.
Otherwise, at least one peak of the dual vector is a sidelobe artifact.  
Then an \emph{additional} peak in $\Vec{u}$ must be included for estimating the next $\mu$. 
This is implemented by incrementing the loop counter $i$ in line 8,
which keeps track of the number of peaks in $\Vec{u}$.}
\delete{
If $(\Mat{A}^H\Mat{A})$ were a diagonal matrix, the relation would also be exact.
This follows from
\begin{eqnarray}
\Vec u^{i} &= & 2 \Mat{D}^{-H} \Mat{A}^H \left(\Vec{y} - \Mat{A} \;  \Vec{x}_{\ell_1}^{i}   \right) \nonumber \\ 
   &=&  2 \Mat{D}^{-H} ( \Mat{A}^H \Vec{y} - (\Mat{A}^H\Mat{A}) \;  \Vec{x}_{\ell_1}^{i}). \label{eq:coherence-effect}
\end{eqnarray}
A diagonal coherence matrix implies the mutual coherence is zero, which is not possible for $M>N$.
From (\ref{eq:coherence-effect}), it is seen that the coherence matrix of the basis $( \Mat{A}^H\Mat{A})$ is re-weighted by the primal solution vector $\Vec{x}_{\ell_1}^{i}$ in the $i$th iteration. The accuracy of the approximation (\ref{eq:fast-approx}) depends on the magnitudes of the off-diagonal elements of the coherence matrix.
For discussing the approximation accuracy, note
that $\mu^i$ depends on  the peaks in the magnitude of $| \Vec u^{i} |$
\begin{eqnarray}
\Vec u^{i}
   &=& 2 \Mat{D}^{-H} \Mat{A}^H \left(\Vec{y} - \sum\limits_{m\in\mathcal{M}_i} \Vec{a}_m \Vec{x}_{\ell_1,m}^{i} \right)
\label{eq:coherence-effect-discussion}
\end{eqnarray}
As long as the \emph{next} active column is orthogonal to 
\emph{all}  active columns $\Vec{a}_m$ with $m\in\mathcal{M}_i$, the approximation is exact. 

Although (\ref{eq:coherence-effect-discussion}) appears to be quite complicated, Corollary 3 assures that the dual coordinate is linear in 
the regularization parameter for all $\mu^{*(p+1)}<\mu<\mu^{*p}$.}
it may happen that the coordinate corresponding to the $(i+1)^{\rm th}$ peak becomes active first,
although $
  \mathop{\mathrm{peak}}(\Vec{u}^{i-1}, i) >   \mathop{\mathrm{peak}}(\Vec{u}^{i-1}, i+1)$.
In this case, two sources become active as the regularization parameter is chosen too low. 
\delete{This is not treated in the tables to keep them simple, but }This exception can be handled by, e.g., bisection in $\mu$.

\begin{table}[h!]
\begin{center}
\begin{tabular}{ll} 
 \hline\hline 
    & Given: $\Mat{A}\in\mathbb{C}^{N\times M}$, $\Mat{D}\in {\rm diag}\mathbb{R}^{M}$, $\Vec{y}\in\mathbb{C}^N$ \\
    & Given: $K_0\in{\mathbb N}$ ,  $F\in]0,1[$, $\Vec{x}_{\ell_1}$. \\ \hline \\
    
1: &   $\mathcal{M} = \{ m \Big|\, |x_{\ell_1,m}| > \delta_i \}$,  $\delta_i= \epsilon \|\Vec{x}_{\ell_1}^i\|_{\infty} $  \\    
2: &   $\Vec{u}^{i-1}=2\Mat{D}^{-H}\Mat{A}^H \left( \Vec{y} - \Mat{A}\Vec{x}_{\ell_1} \right)  $ \\ \\
3: &   \text{if} $|\mathcal{M}|<K_0$ \\
4:    &   \hspace{2ex} $\mathcal{U} = \{ m \Big|\, 1-\frac{|u_m|}{\mu}<\epsilon_{\mu}\}$\\
5: &  \hspace{2ex} $i = |\mathcal{U}|+1$   \\
6: & \hspace{2ex} $\mu =  (1-F)\mathop{\mathrm{peak}}\! \left(\Vec{u}^{i-1} , i \right)+ F\mathop{\mathrm{peak}}\! \left(\Vec{u}^{i-1} , i+1 \right) $\\
7: & \text{else if} $|\mathcal{M}|>K_0$ \\
8: &  \hspace{2ex} bisecting	 between $\mu^{i-1}$ and $\mu^{i}$ as defined in Eq. (\ref{eq:fix_point}) \\ 
9: & \text{end}\\ \\
10: & Output: $\mu$
\\ \\
\hline\hline
\end{tabular}
\caption{Order-recursive algorithm to select $\mu$ for given sparsity order $K_0$. } 
\label{algo:order-recursive_new}
\end{center}
\end{table}

The recursive-in-order algorithm provided in Table \ref{algo:order-recursive_new} takes as input the dictionary $\Mat{A}$, the generalization matrix $\Mat{D}$, the measurement vector $\Vec{y}$, the given sparsity order $K_0$ and the previous order LASSO solution $\Vec{x}_{\ell_1}$. In line 1 the actual active set is determined by thresholding and line 2 produces the dual vector by Theorem \ref{theo:dual}. Line 2 can be omitted, if the LASSO solver makes the dual solution available, e.g., through primal-dual interior point methods or alternating direction method of multipliers. If the size of the active set of the previous LASSO solution is less than the given sparsity order $K_0$, the algorithm determines the dual active set $\mathcal{U}$ in line 4, cf. Eq.(\ref{eq:active-set}). The incremented cardinality of $\mathcal{U}$ is the new requested number of hitting peaks in the dual vector. Finally, line 6 calculates $\mu$ based on the candidate point estimate  (\ref{eq:mu-approx}).

\subsection{Fast-Iterative Algorithm}
\label{se:fast}

The approximation from Equation (\ref{eq:mu-approx}) is not limited to a single iteration. 
Therefore, (\ref{eq:mu-approx}) can be extended further to
\begin{align}
\mu^i  &\approx {\rm peak}(\Vec{u}(\mu^{i-1}), i)  \notag
                 \\ &\approx {\rm peak}(\Vec{u}(\mu^{i-2}), i)  \notag \\ & \approx
         \cdots \notag
\\ & \approx  {\rm peak}(\Vec{u}(\mu^{0}), i)= {\rm peak}(2 \Mat{D}^{-H} \Mat{A}^H \Vec{y}, i)~. \label{eq:fast-approx}
\end{align}
This  observation motivates the iterative algorithm  in Table \ref{algo:fast_new}. The main difference to the recursive-in-order algorithm is found in line 6. The peakfinder estimates the maximum of the $K^{\rm th}$ peak. This leads to a significant speed-up, if sources are well separated and their basis coherence is low.

\begin{table}[h!]
\begin{center}
\begin{tabular}{ll} 
 \hline\hline 
    & Given: $\Mat{A}\in\mathbb{C}^{N\times M}$, $\Mat{D}\in {\rm diag}\mathbb{R}^{M}$, $\Vec{y}\in\mathbb{C}^N$ \\
    & Given: $K_0\in{\mathbb N}$ ,  $F\in]0,1[$, $\Vec{x}_{\ell_1}$ . \\ \hline \\
    
1: &   $\mathcal{M} = \{ m \Big|\, |x_{\ell_1,m}| > \delta_i \}$,  $\delta_i= \epsilon \|\Vec{x}_{\ell_1}^i\|_{\infty} $  \\    
2: &   $\Vec{u}^{i-1}=2\Mat{D}^{-H}\Mat{A}^H \left( \Vec{y} - \Mat{A}\Vec{x}_{\ell_1} \right)  $ \\ \\
3: &   \text{if} $|\mathcal{M}|<K_0$ \\
4:    &   \hspace{2ex} $\mathcal{U} = \{ m \Big|\, 1-\frac{|u_m|}{\mu}<\epsilon_{\mu}\}$\\
5: &  \hspace{2ex} $i = |\mathcal{U}|+1$   \\
6: & \hspace{2ex} $\mu =  (1-F)\mathop{\mathrm{peak}}\! \left(\Vec{u}^{i-1} , K \right)+ F\mathop{\mathrm{peak}}\! \left(\Vec{u}^{i-1} , K_0+1 \right) $\\
7: & \text{else if} $|\mathcal{M}|>K_0$ \\
8: &  \hspace{2ex} bisecting	 between $\mu^{i-1}$ and $\mu^{i}$ as defined in Eq. (\ref{eq:fix_point})		\\
9: & \text{end}\\ \\
10: & Output: $\mu$
\\ \\
\hline\hline
\end{tabular}
\caption{Iterative primal based algorithm to select $\mu$ for given sparsity order $K_0$. }
\label{algo:fast_new}
\end{center}
\end{table}

\subsection{Detection in the dual domain}

As a demonstrative example, we provide the fast iterative algorithm formulated \emph{solely} in the dual domain in Table \ref{algo:dual-fast_new}. Note that the gird-free atomic norm solutions \cite{Tang2013,Candes2013,Candes2014,Panahi2014icassp,Xenaki2015} follow a similar approach.

As asserted by (\ref{eq:active-set}), searching for active indices in the dual domain is effectively a form of relaxation of the primal problem (\ref{eq:lasso_prob}). This amounts to peak finding in the output of a WMF acting on the residuals, cf. Theorem \ref{theo:dual}.
In line 1, the active set $\mathcal{M}$ is effectively approximated by the relaxed set $\mathcal{U}$. Therefore, the $\ell_0$ solution is determined by regression on the relaxed set in line 2 and the primal active set is found by thresholding this solution in line 3. The remainder of the algorithm is equal the primal based ones.

\delete{Instead of the primal (\ref{eq:lasso_prob}),  the dual  (\ref{eq:tibshirani-dual-problem35-min})--(\ref{eq:row-space-constraint})  is solved exclusively.}

\begin{table}[h!]
\begin{center}
\begin{tabular}{ll} 
 \hline\hline 
    & Given: $\Mat{A}\in\mathbb{C}^{N\times M}$, $\Mat{D}\in {\rm diag}\mathbb{R}^{M}$, $\Vec{y}\in\mathbb{C}^N$ \\
    & Given: $K_0\in{\mathbb N}$ ,  $F\in]0,1[$, $\Vec{u}$. \\ \hline \\

1: &  $\mathcal{U} = \{ m \Big|\, 1-\frac{|u_m|}{\mu}<\epsilon_{\mu}\}$ \\ 
2: & $\Vec{x}_{\ell_0} = \Mat{A}_{\mathcal{U}}^+ \Vec{y}$ \\ 
3: &    $\mathcal{M} = \{ m \Big|\, |x_{\ell_0,m}| > \delta \}$,  $\delta= \epsilon \|\Vec{x}_{\ell_0}\|_{\infty} $  \\ \\

4: &   \text{if} $|\mathcal{M}|<K_0$ \\
5: &  \hspace{2ex} $i = |\mathcal{U}|+1$   \\
6: & \hspace{2ex} $\mu =  (1-F)\mathop{\mathrm{peak}}\! \left(\Vec{u}^{i-1} , K \right)+ F\mathop{\mathrm{peak}}\! \left(\Vec{u}^{i-1} , K_0+1 \right) $\\
7: & \text{else if} $|\mathcal{M}|>K_0$ \\
8: &  \hspace{2ex} bisecting	 between $\mu^{i-1}$ and $\mu^{i}$ as defined in Eq. (\ref{eq:fix_point})			\\
9: & \text{end}\\ \\
10: & Output: $\mu$
\\ \\
\hline\hline
\end{tabular}
\caption{Iterative dual based algorithm to select $\mu$ for given sparsity order $K_0$.}
\label{algo:dual-fast_new}
\end{center}
\end{table}

\section{Simulation}
In this section, the performance of the proposed dual estimation algorithms is evaluated based on numerical simulation.
We use synthetic data from a uniform linear array with $N=64$ elements with half-wavelength spacing. 
The DOA domain is discretized by  $\theta_m = (m-1)\frac{180^\circ}{M}-90^{\circ}$ with $m=1,\ldots,M$ and $M=180$.
The simulation scenario has $K_0=8$ far-field  plane-waves sources (\ref{eq:linear-model2}). 
The  uncorrelated noise $\Vec{n}$ is  zero-mean complex-valued circularly symmetric normally distributed $\sim\mathcal{N}(\Vec{0},\,\Vec{I})$, i.e. $0\,$dB power.
Eight sources are stationary at 
$\theta^{\rm T}= [-45,\,   - 30,\,   -14,\,    9,\,   17,\,   30,\,   44,\, 72]$ 
degrees relative to broadside with constant power level (PL) $[-5,\, 10,\, 5,\, 0,\, 11,\, 12,\, 9,\, 25]$\,dB \cite{MGPV2013}.

\begin{figure} 
\centering
	\includegraphics[width=1\columnwidth]{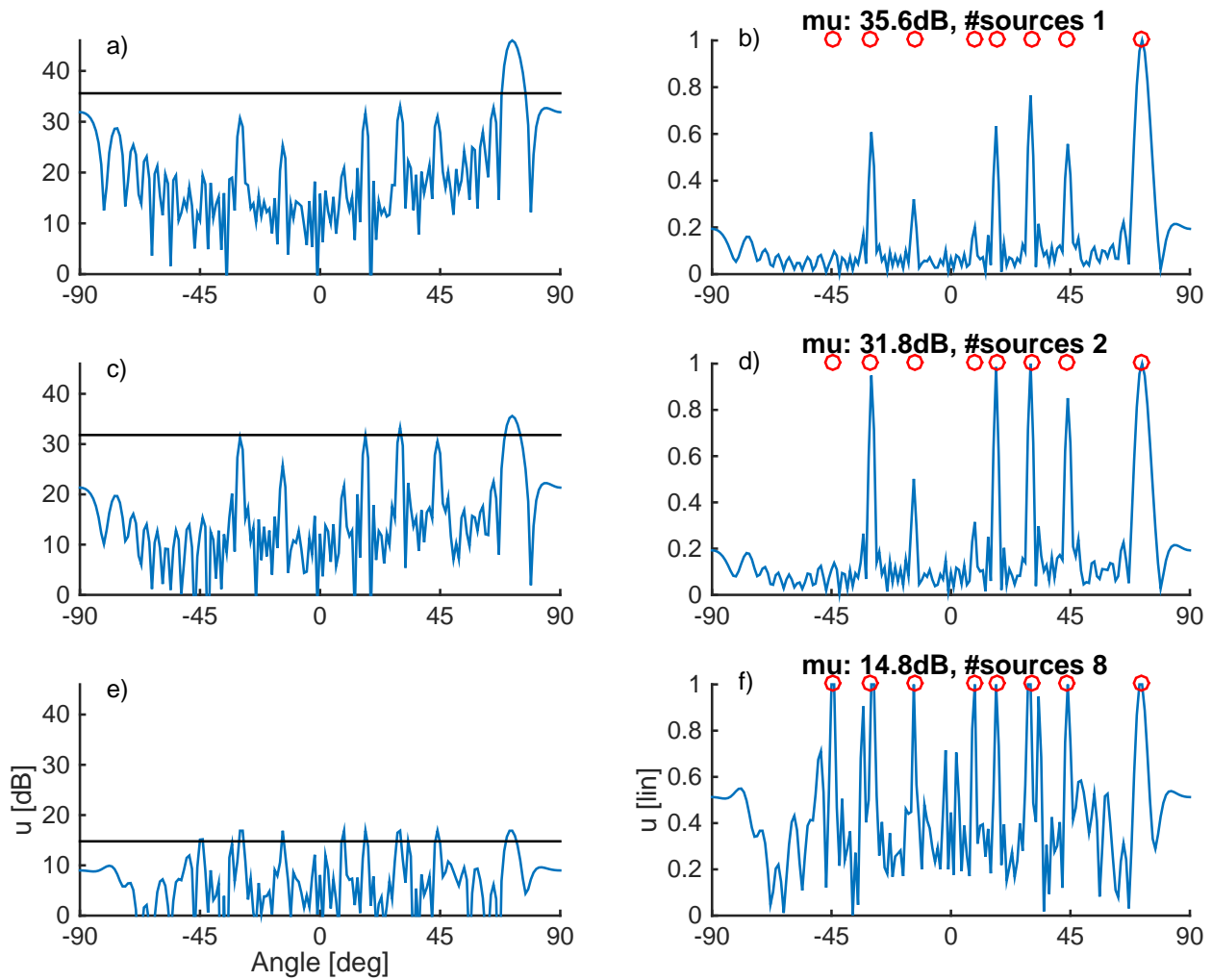}
\caption{Dual  coordinates  for order-recursive approach corresponding to step $p=1$ (a and b),  $p=2$ (c and d), and  $p=8$ (e and f).
Left column: Dual (dB) for the previous step which is used for selecting $\mu$ (horizontal line).
Right column: Dual (lin) normalized  with $\mu$ (maximum is 1), the true source locations are marked with $\circ$, and the actual value of $\mu$ and number of sources found is also indicated.
} \label{fig:dualmethod}
\end{figure}

\begin{figure} 
\centering
	\includegraphics[width=1\columnwidth]{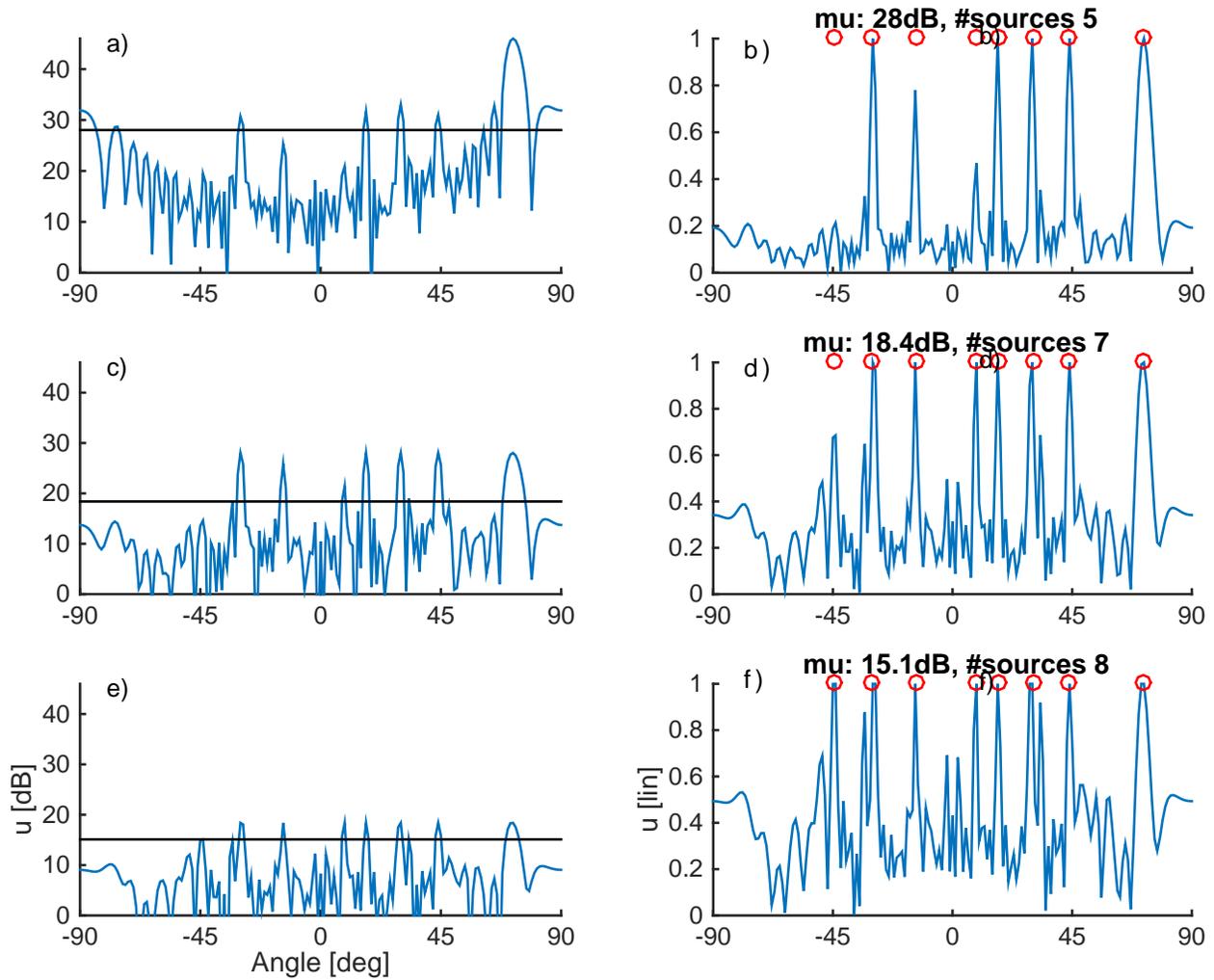}
\begin{picture}(0,0)
\put(19,134){\small\sf b\,)}
\put(19,90){\small\sf d\,)}
\put(19,47){\small\sf f\,)}
\end{picture}
\caption{Dual  coordinates  iterative approach corresponding for localizing $K_0=8$ sources for step $i=1$ (a and b), $i=2$ (c and d), and $i=3$ (e and f).
Left column: Dual (dB) for the previous step which is used for selecting $\mu$ (horizontal line).
Right column: Dual (lin) normalized  with $\mu$ (maximum is 1), the true source locations are marked with $\circ$, and the actual value of $\mu$ and number of sources found is also indicated..
} \label{fig:dualmethoditer}
\end{figure}

The dual solution for the order-recursive approach, Table \ref{algo:order-recursive_new}, corresponds to the results shown in Fig.\  \ref{fig:dualmethod}. 
The faster iterative approach, Table \ref{algo:fast_new}, yields the results in Fig.\  \ref{fig:dualmethoditer}. 
The dual solution using the primal solution from the previous iteration is interpreted as a WMF and used for the selection of $\mu$ (left column). 
Next, the convex optimization is carried out for that value of $\mu$ giving the dual solution. 
We plot the dual solution on a linear scale and normalized to a maximum value of 1 which is customary in implementations of the dual for compressed sensing \cite{Tang2013,Candes2013,Candes2014}. 
The number of active sources (see right column in Figs.\  \ref{fig:dualmethod}  and \ref{fig:dualmethoditer}) 
are determined according to line 1 in Tables \ref{algo:order-recursive_new} and \ref{algo:fast_new}.

For the order-recursive approach step 1, Fig.\  \ref{fig:dualmethod}a, the $\mu $ is selected based on the main peak $\theta=72^\circ$ and a large side lobe at $\theta=80^{\circ}$. Once the solution for that $\mu$ is obtained it turns out that there   is no  an active source in the sidelobe.The solution progresses steadily down the LASSO path.
Figure  \ref{fig:dualmethoditer} shows the faster iterative approach in Table \ref{algo:fast_new} for the 8-source problem. In the first iteration we use a $\mu$ between the 8th and 9th peak based on the WMF solution (Fig.\ \ref{fig:dualmethoditer}a). There are many sidelobes associated with the source at $\theta=72^\circ$. 
As soon as the dominant source is determined, the sidelobes in the residuals are reduced and only 5 sources are observed. 
After two more iterations, all 8 sources are found at their correct locations.

For both algorithms, the main CPU time is used in solving the convex optimization problem. 
Thus the iterative algorithm is a factor 8/3 faster in this case than the straightforward approach which strictly follows the LASSO path.
The approach described in Table \ref{algo:order-recursive_new} has approximately the same CPU time usage as  the approach in Ref. \cite{MGPV2013}, but it is conceptually simpler 
and provides deeper physical insight into the problem.


\section{Conclusion}
The complex-valued generalized  LASSO problem is convex. 
The corresponding dual problem is  interpretable as a weighted matched Filter (WMF) acting on the residuals of the LASSO.
There is a linear one-to-one relation between the dual and primal vectors. Any results formulated for the primal problem are readily extendable to the dual problem.
Thus, the sensitivity of the primal solution to small changes in the constraints can be easily assessed.
Further, the difference between the solutions $\Vec{x}_{\ell_0}$ and the $\Vec{x}_{\ell_1}$ is characterized via the dual vector.

Based on mathematical and physical insight,  an order-recursive and a faster iterative LASSO-based algorithm are proposed and evaluated.
These algorithms use the dual variable of the generalized LASSO for regularization parameter selection. This greatly facilitates computation of the LASSO-path as we can predict the changes in the active indexes as the regularization parameter is reduced.
Further, a dual-based algorithm is formulated which solves only the dual problem.
The examples demonstrate the algorithms, confirming that the dual and primal coordinates are piecewise linear in the regularization parameter $\mu$.

\section*{Appendix A}
\renewcommand{\theequation}{A\arabic{equation}}
\setcounter{equation}{0}

Proof of (\ref{eq:auxiliary-lemma}): Set $\Vec{u}=(u_1,\ldots,u_M)^T\in\mathbb{C}^M$. 
From (\ref{eq:L2-involving-z}),
\begin{eqnarray}\label{eq:fortyseven} \hspace{-2ex}
  \mu \| \Vec{z} \|_1 - \Re(\Vec{u}^H\Vec{z}) \hspace{-1.5ex} &=& \hspace{-2ex}
 \sum\limits_{m=1}^M  \left( \mu |z_m| - \Re(u_m^*z_m) \right) \\
 \hspace{-1.5ex} &=& \hspace{-2ex} \sum\limits_{m=1}^M \underbrace{(\mu  - |u_m|\cos\phi_{mm})}_{=\tilde{\mu}_m} |z_m|, \label{eq:fortyseven2}
\end{eqnarray}
where we set $u_m^*z_m=|u_m|\,|z_m|\,\mathrm{e}^{j\phi_{mm}}$. 
The phase difference $\phi_{mm}$ depends on both $u_m$ and $z_m$. 
If all coefficients  $\tilde{\mu}_m$ in (\ref{eq:fortyseven2}) are non-negative, $\tilde{\mu}_m \ge 0$, for all $z_m\in\mathbb{C}$, then
\be
    \min\limits_{\Vec{z}} \left( \mu \| \Vec{z} \|_1 - \Re(\Vec{u}^H\Vec{z})\right) = 0,
\ee
otherwise there is no lower bound on the minimum. Therefore, all $|u_m|$ must be bounded, i.e. $ |u_m| \le \mu\,\forall \, m=1,\ldots,M$ 
to ensure that all $\tilde{\mu}_m\ge0$ for all possible phase differences $-1\le\cos\phi_{mm}\le1$.
Finally, we note that 
$\|\Vec{u}\|_{\infty}=\max_m|u_m|$.

\section*{Appendix B: Proofs of Corollaries \ref{cor:boundary}, \ref{cor:phase}, and \ref{cor:piecewise-linear}}
\renewcommand{\theequation}{B\arabic{equation}}
\setcounter{equation}{0}

\subsection*{Proof for Corollary \ref{cor:boundary}}
Let the objective function of the complex-valued generalized LASSO problem (\ref{eq:lasso_prob}) be
\be
     \mathscr{L} = \| \Vec{y} - \Mat{A} \Vec{x} \|_2^2 + \mu \| \Mat{D}\Vec{x} \|_1~.
\ee
In the following, we evaluate the subderivative $\partial\mathscr{L}$ \cite{Bertsekas1999} as the set of all complex subgradients as introduced in \cite{Bouboulis2012}. First, we observe
\be\label{eq:subdifferential}
    \partial\mathscr{L} = -2\Mat{A}^H(\Vec{y}-\Mat{A}\Vec{x}) + \mu \,\partial\| \Mat{D} \Vec{x} \|_1~.
\ee
Next, it is assumed that $\Mat{D}$ is a diagonal matrix with positive real-valued diagonal entries. Then the subderivate $\partial \| \Mat{D}\Vec{x} \|_1$ evaluates to
\be 
\partial \| \Mat{D}\Vec{x} \|_1 = \left\{\begin{array}{cc}
\frac{D_{mm} x_m}{|x_m|} & \mathrm{for } \quad x_{m} \neq 0 \\ 
& \\
\{z \in \mathbb{C}, |z|\leq 1 \}& \mathrm{for } \quad x_m = 0.
\end{array}   \right.
\ee
The minimality condition for $\mathscr{L}$ is equivalent to setting (\ref{eq:subdifferential}) to zero. For all  $m$ with $ x_{m} \neq 0 $
and with (\ref{eq:2nd_time_mystery}), this gives
\be \label{eq:subgrad_active_set}
D_{mm} u_m = \mu \frac{D_{mm}x_m}{|x_m|}.
\ee 
It readily follows that $|u_m|=\mu$ for $x_m\ne0$ and $D_{mm}\ne0$.


\subsection*{Proof for Corollary \ref{cor:phase}}
Starting from Eq. (\ref{eq:subgrad_active_set}), dividing by $\mu$ and invoking Corollary 1, we conclude for matrices $\Mat{D}$ with positive diagonal entries
and for $m\in\mathcal{M}$,
\be \label{eq:subgrad_active_set2}
\mu \mathrm{e}^{j\arg(x_m)} = \frac{2}{D_{mm}} \Vec{e}_m^H\Mat{A}^H \left( \Vec{y} - \Mat{A}\Vec{x}\right) = u_m ~,
\ee
where $\Vec{e}_m$ is the $m$th standard basis vector.
This concludes the proof of Corollary \ref{cor:phase}.

\subsection*{Proof for Corollary \ref{cor:piecewise-linear}}
For the primal vector, this was shown in the real-valued case by Tibshirani \cite{Tibshirani1996} and for the complex-valued case, this is a direct consequence of Appendix B in \cite{MGPV2013}.
For the dual vector, this was shown in the real-valued case by Tibshirani \cite{Tibshirani2011} and for the complex-valued case, this readily follows from Theorem 1: If the primal vector $\Vec{x}_{\ell_1}$ depends linearly on $\mu$ in (\ref{eq:2nd_time_mystery}) then so does the dual vector $\Vec{u}$.

\section*{Appendix C: $\ell_0$ solution}
\renewcommand{\theequation}{C\arabic{equation}}
\setcounter{equation}{0}

The gradient (cf. Appendix B) of the 
data objective function is
\be
\nabla \left\| \Vec{y} - \Mat{A} \Vec{x}  \right\|_2^2
= -2 \Mat{A}^H \left( \Vec{y} - \Mat{A} \Vec{x} \right) \label{eq:grad_almost_u}
\ee
For the active source components, $x_m$ with $m\in\mathcal{M}$, 
the $\ell_0$-constraint of (\ref{eq:l0_prob}) is without effect and the solution results from setting the gradient to zero, i.e. solving the normal equations. 
\begin{align}
 \Mat{A}_\mathcal{M}^H  \Vec{y} =  \Mat{A}_\mathcal{M}^H\Mat{A}_\mathcal{M} \Vec{x}_{\ell_0,\mathcal{M}}  \quad\Rightarrow \quad\Vec{x}_{\ell_0,\mathcal{M}}=\Mat{A}_\mathcal{M}^+ \Vec{y}
\end{align}
\\
We set
\be
\Vec{x}_{\ell_0,\mathcal{M}} = \Vec{x}_{\ell_1,\mathcal{M}}+\Vec{\Delta} ~.
\label{eq:Delta}
\ee
This is inserted into (\ref{eq:grad_almost_u}),
\begin{eqnarray}
\nabla \left\| \Vec{y} - \Mat{A} \Vec{x}_{\ell_1,\mathcal{M}}  \right\|_2^2
= -2 \Mat{A}^H \left( \Vec{y} - \Mat{A} (\Vec{x}_{\ell_0,\mathcal{M}} - \Vec{\Delta}) \right).
\end{eqnarray}
Using (\ref{eq:Tibshirani-Eq34-Erich}) gives
\begin{eqnarray}
 \Mat{D}_\mathcal{M}^{H} \Vec{u}_\mathcal{M}   &=&   2 \Mat{A}_\mathcal{M}^H \left( \Vec{y} - \Mat{A}_\mathcal{M} \Vec{x}_{\ell_1,\mathcal{M}} \right)   \\
\Mat{D}_\mathcal{M}^{H} \mu e^{j\Vec{\theta}_\mathcal{M}}     & = &       2 \Mat{A}_\mathcal{M}^H \left( \Vec{y} - \Mat{A}_\mathcal{M} \left(\Vec{x}_{\ell_0,\mathcal{M}} - \Vec{\Delta}\right) \right) \\
 \mu  \Mat{D}_\mathcal{M}^{H} e^{j\Vec{\theta}_\mathcal{M}}     & =&       2 \Mat{A}_\mathcal{M}^H \Mat{A}_\mathcal{M} \Vec{\Delta}  
\end{eqnarray}
This results in
\be
\Vec{\Delta} = \frac{\mu}{2} \left(\Mat{A}_\mathcal{M}^H \Mat{A}_\mathcal{M}\right)^{-1}\Mat{D}_\mathcal{M}^H
 e^{j\Vec{\theta}_\mathcal{M}} \label{eq:Delta-result}
\ee
which depends on $\mu$ both explicitly and implicitly through $\mathcal{M}$.
If the set of nonzero elements of (\ref{eq:l0_prob}) is equal to the active set of (\ref{eq:lasso_prob}), the solutions of (\ref{eq:l0_prob}) and (\ref{eq:lasso_prob})
differ by (\ref{eq:Delta-result}).


\end{document}